\documentclass[a4paper,12pt]{article}

\usepackage{latexsym}
\usepackage{amssymb}
\usepackage{theorem}
\usepackage{amsmath}
\usepackage{amscd}
\usepackage{graphicx}
\usepackage{color}

\newtheorem{theorem}{Theorem}[section]

\newtheorem{lemma}[theorem]{Lemma}
\newtheorem{example}[theorem]{Example}
\newtheorem{proposition}[theorem]{Proposition}
\newtheorem{remark}[theorem]{Remark}
\newtheorem{definition}[theorem]{Definition}

\newcommand{\demo}{\par\noindent{\it Proof. \/}\ }
\newcommand{\enD}{\hfill $\Box$\vspace{3truemm} \par}
\newcommand{\R}{\mathbb{R}}

\newcommand{\bn}{\mbox{\boldmath $n$}}
\newcommand{\bt}{\mbox{\boldmath $t$}}
\newcommand{\ba}{\mbox{\boldmath $a$}}
\newcommand{\bb}{\mbox{\boldmath $b$}}
\newcommand{\be}{\mbox{\boldmath $e$}}

\newcommand{\bv}{\mbox{\boldmath $v$}}
\newcommand{\bw}{\mbox{\boldmath $w$}}
\newcommand{\bx}{\mbox{\boldmath $x$}}

\newcommand{\bxi}{\mbox{\boldmath $\xi$}}
\newcommand{\beeta}{\mbox{\boldmath $\eta$}}
\newcommand{\Ev}{\mathcal{E}v}
\newcommand{\Inv}{\mathcal{I}nv}

\begin{document}

\title{Evolutes and involutes of framed curves in Euclidean 3-space}

\author{Nozomi Nakatsuyama}

\date{\today}

\maketitle

\begin{abstract}
We investigated the evolute of a space curve with singular points. As smooth curves with singular points, we apply the theory of framed curves.
However, the involute corresponding to the evolute in the sense of the locus of the centre of osculating spheres has not been defined as far as we know. 
In this paper, we directly define the evolutes and involutes of non-degenerate curves and framed curves using the theory of Bertrand type curves. 
We give conditions that the evolutes and involutes are inverse operations of these curves. Moreover, we investigate the other cases in which the evolutes and involutes are not inverse operations. 
\end{abstract}

\renewcommand{\thefootnote}{\fnsymbol{footnote}}

\section{Introduction}
The notions of evolutes and involutes (also known as evolvents) are classical objects in differential geometry.
There are many articles concerning the evolute (i.e. the locus of the centre of osculating spheres) of a Frenet curve (i.e. a regular space curve with the linear independent condition) for instance, \cite{Arnold1, Blaschke, Fuchs, Gray, Porteous, Romero-Sanabria, Uribe-Vargas1}. 
In \cite{Honda-Takahashi-2020-evolute}, they investigated the evolute of a space curve with singular points. 
As smooth curves with singular points, we apply the theory of framed curves (cf. \cite{Honda-Takahashi-2016}). 
A framed curve is a space curve with a moving frame. 
Framed curves may have singular points, and have a circular evolute and an evolute in the sense of the locus of the centre of osculating spheres. 
The circular evolutes and involutes of framed curves have already been investigated (cf. \cite{Honda-Takahashi-Preprint}). 
However, the involute corresponding to the evolute has not been defined as far as we know. 
\par
On the other hand, Bertrand and Mannheim curves are also classical objects in differential geometry (cf. \cite{Aminov, Banchoff-Lovett, Berger-Gostiaux, Bertrand, doCarmo, HCIP, Izumiya-Takeuchi1, Kuhnel, Lucas-Ortega, Struik}). 
A Bertrand curve (respectively, Mannheim curve) is a space curve whose principal normal line is the same as the principal normal line (respectively, the bi-normal line) of another curve.
In \cite{Nakatsuyama-Takahashi}, we investigated the existence conditions of $(\bv,\overline{\bw})$-Bertrand type curves in all cases. 
As a consequence, the planar evolutes and planar involutes appear as the $(\bv,\overline{\bw})$-Bertrand type curves.
On the other hand, we defined and investigated the existence conditions of the $(\bv,\overline{\bw})$-Bertrand framed curves in all cases. 
As a consequence, the circular evolutes and involutes of framed curves appear as the $(\bv,\overline{\bw})$-Bertrand framed curves. 
It follows that the evolutes and involutes of a non-degenerate curve (respectively, a framed curve) are considered to be found in Bertrand type curves (respectively, Bertrand framed curves). 
However, evolutes in the sense of the locus of the centre of osculating spheres does not appear. 
\par
In this paper, we generalize the notions of Bertrand type curves and Bertrand framed curves.
Then, we define the evolutes and involutes of non-degenerate curves (respectively, framed curves) by using the method of $(\bv,\bw,\overline{\bx})$-Bertrand type curves (respectively, $(\bv,\bw,\overline{\bx})$-Bertrand framed curves). 
In \S 2, we review the theories of Bertrand and Mannheim curves of non-degenerate curves and framed curves (cf. \cite{Honda-Takahashi-2020}). 
We define and clarify that the existence conditions of the Bertrand type curves and Bertrand framed curves in \S 3 and \S 4. 
As a consequence, the evolutes and involutes of non-degenerate curves (cf. \cite{Blaschke, Bruce-Giblin, Gibson,Gray}) appear as one of Bertrand type curves (Theorems \ref{nbb-Bertrand-type} and \ref{tnt-Bertrand-type}). 
Moreover, the evolutes and involutes of framed curves (cf. \cite{Honda-Takahashi-2020-evolute}) appear as one of Bertrand framed curves (Theorems \ref{n1n2n2-Bertrand-type} and \ref{mun1mu-Bertrand-type}). 
In \S 5 and \S 6, we directly define the evolutes and involutes of non-degenerate curves (respectively, framed curves) using the theory of Bertrand type curves (respectively, Bertrand framed curves). 
We give the relation between the evolute (respectively, involute) of a framed curve and the circular evolute with respect to Bishop directions of a framed curve (respectively, involute of $\gamma$ respect to $t_0$) (Proposition \ref{evolute-evolute} (respectively, Proposition \ref{involute-involute})). 
We give conditions that the evolutes and involutes are inverse operations of these curves (Theorems \ref{evo_invo_non-deg_inverse} and \ref{evo_invo_framed_inverse}). Moreover, we investigate other cases in which the evolutes and involutes are not inverse operations (Theorems \ref{Inv(Ev)_another} and \ref{evo_invo_framed_inverse-another}). 
\par
We shall assume throughout the whole paper that all maps and manifolds are $C^{\infty}$ unless the contrary is explicitly stated.

\bigskip
\noindent
{\bf Acknowledgement}. 
The author would like to express sincere gratitude to Professor Kentaro Saji and Professor Masatomo Takahashi for helpful comments. This work was supported by JST SPRING Grant Number JPMJSP2153.

\section{Preliminaries}

We review the theories of Bertrand and Mannheim of non-degenerate curves and framed curves.

Let $\R^3$ be the $3$-dimensional Euclidean space equipped with the inner product $\ba \cdot \bb = a_1 b_1 + a_2 b_2 + a_3 b_3$, 
where $\ba = (a_1, a_2, a_3)$ and $\bb = (b_1, b_2, b_3) \in \R^3$. 
The norm of $\ba$ is given by $\vert \ba \vert = \sqrt{\ba \cdot \ba}$ and the vector product is given by 
$$
\ba \times \bb={\rm det}
\left(
\begin{array}{ccc}
\be_1 & \be_2 & \be_3\\
a_1 & a_2 & a_3\\
b_1 & b_2 & b_3
\end{array}
\right),
$$
where $\{\be_1, \be_2, \be_3\}$ is the canonical basis of $\R^3$. 
Let $S^2$ be the unit sphere in $\R^3$, that is, $S^2=\{\ba \in \R^3| |\ba|=1\}$.
We denote the $3$-dimensional smooth manifold $\{(\ba,\bb) \in S^2 \times S^2| \ba \cdot \bb=0\}$ by $\Delta$. 

\subsection{Regular curves}

Let $I$ be an interval of $\R$ and let $\gamma:I \to \R^3$ be a regular space curve, that is, $\dot{\gamma}(t) \not=0$ for all $t \in I$, where $\dot{\gamma}(t)=(d\gamma/dt)(t)$. 
We say that $\gamma$ is {\it non-degenerate}, or $\gamma$ satisfies the {\it non-degenerate condition} if $\dot{\gamma}(t) \times \ddot{\gamma}(t) \not=0$ for all $t \in I$. 

If we take the arc-length parameter $s$, that is, $|\gamma'(s)|=1$ for all $s$, then the tangent vector, the principal normal vector and the bi-normal vector are given by
$
\bt(s)=\gamma'(s), \ \bn(s)={\gamma''(s)}/{|\gamma''(s)|}, \ \bb(s)=\bt(s) \times \bn(s),
$
where $\gamma'(s)=(d\gamma/ds)(s)$. 
Then $\{\bt(s),\bn(s),\bb(s)\}$ is a moving frame of $\gamma(s)$ and we have the Frenet-Serret formula: 
$$
\left(
\begin{array}{c}
\bt'(s)\\
\bn'(s)\\
\bb'(s)
\end{array}
\right)
=
\left(
\begin{array}{ccc}
0&\kappa(s)&0\\
-\kappa(s)&0&\tau(s)\\
0&-\tau(s)&0
\end{array}
\right)
\left(
\begin{array}{c}
\bt(s)\\
\bn(s)\\
\bb(s)
\end{array}
\right),
$$
where 
$\kappa(s)=|\gamma''(s)|, \ \tau(s)={{\rm det}(\gamma'(s),\gamma''(s),\gamma'''(s))}/{\kappa^2(s)}$.
If we take a general parameter $t$, then the tangent vector, the principal normal vector and the bi-normal vector are given by
$
\bt(t)={\dot{\gamma}(t)}/{|\dot{\gamma}(t)|}, \ \bn(t)=\bb(t) \times \bt(t), \ \bb(t)={(\dot{\gamma}(t) \times \ddot{\gamma}(t))}/{|\dot{\gamma}(t) \times \ddot{\gamma}(t)|}.
$
Then $\{\bt(t),\bn(t),\bb(t)\}$ is a moving frame of $\gamma(t)$ and we have the Frenet-Serret formula: 
$$
\left(
\begin{array}{c}
\dot{\bt}(t)\\
\dot{\bn}(t)\\
\dot{\bb}(t)
\end{array}
\right)
=
\left(
\begin{array}{ccc}
0&|\dot{\gamma}(t)|\kappa(t)&0\\
-|\dot{\gamma}(t)|\kappa(t)&0&|\dot{\gamma}(t)|\tau(t)\\
0&-|\dot{\gamma}(t)|\tau(t)&0
\end{array}
\right)
\left(
\begin{array}{c}
\bt(t)\\
\bn(t)\\
\bb(t)
\end{array}
\right),
$$
where 
$
\kappa(t)={|\dot{\gamma}(t) \times \ddot{\gamma}(t)|}/{|\dot{\gamma}(t)|^3}, \ \tau(t)={{\rm det}(\dot{\gamma}(t),\ddot{\gamma}(t),\dddot{\gamma}(t))}/{|\dot{\gamma}(t) \times \ddot{\gamma}(t)|^2}
$.
Note that in order to define $\bn(t), \bb(t)$ and $\tau(t)$, we assume that $\gamma$ is not only regular, but also non-degenerate.
\begin{definition}\label{non-deg_evolute}
The evolute $Ev(\gamma):I \rightarrow \mathbb{R}^3$ of a non-degenerate curve $\gamma$ with $\tau\neq0$ is given by 
$$
Ev(\gamma)(t) = \gamma(t) + \frac{1}{\kappa(t)} \bn(t) - \frac{\dot{\kappa}(t)}{\vert \dot{\gamma}(t) \vert \kappa^2(t) \tau(t)} \bb(t). 
$$
\end{definition}
\subsection{Bertrand and Mannheim non-degenerate curves}

Let $\gamma$ and $\overline{\gamma}:I \to \R^3$ be different non-degenerate curves. 
\begin{definition}\label{regular-Bertrand-Mannheim}{\rm 
We say that $\gamma$ and $\overline{\gamma}$ are {\it Bertrand mates} (respectively, {\it Mannheim mates}) if there exists a smooth function $\lambda:I \to \R$ such that $\overline{\gamma}(t)=\gamma(t)+\lambda(t)\bn(t)$ and $\bn(t)=\pm \overline{\bn}(t)$ (respectively, $\bn(t)=\pm \overline{\bb}(t)$) for all $t \in I$.
\par
We also say that $\gamma:I \to \R^3$ is a {\it Bertrand curve} (respectively, {\it Mannheim curve}) if there exists another non-degenerate curve $\overline{\gamma}:I \to \R^3$ such that $\gamma$ and $\overline{\gamma}$ are Bertrand mates  (respectively, Mannheim mates).}
\end{definition}
By a parameter change, we may assume that $s$ is the arc-length parameter of $\gamma$.
The proof of the existence conditions of Bertrand and Mannheim curves of non-degenerate curves see \cite{Honda-Takahashi-2020} and \cite{Nakatsuyama-Takahashi}.
\begin{theorem}[\cite{Nakatsuyama-Takahashi}]\label{regular-Bertrand-Mannheim-condition}
Let $\gamma:I \to \R^3$ be non-degenerate with the arc-length parameter. 
\par
$(1)$ Suppose that there exists a point $s_0 \in I$ such that $\tau(s_0) \not=0$. 
Then $\gamma$ is a Bertrand curve if and only if there exists a non-zero constant $A$ and a constant $B$ such that  $A\kappa(s)+B\tau(s)=1$ and $\tau(s) (B \kappa(s)-A \tau(s)) \not=0$ for all $s \in I$. 
\par
$(2)$ Suppose that $\tau(s)=0$ for all $s \in I$. 
Then $\gamma$ is a Bertrand curve. 
\par 
$(3)$ $\gamma$ is a Mannheim curves if and only if there exists a non-zero constant $A$ such that  $A(\kappa^2(s)+\tau^2(s))=\kappa(s)$ and $\tau(s)(\kappa(s)\tau'(s)-\kappa'(s)\tau(s)) \not=0$ for all $s \in I$.
\end{theorem}

\subsection{$(\bv,\overline{\bw})$-Bertrand types of non-degenerate curves}
Let $\gamma$ and $\overline{\gamma}:I \to \R^3$ be non-degenerate curves. 
\begin{definition}\label{Bertrand-type-regular}{\rm
We say that $\gamma$ and $\overline{\gamma}$ are {\it $(\bv,\overline{\bw})$-mates} if there exists a smooth function $\lambda:I \to \R$ with $\lambda \not\equiv 0$ such that $\overline{\gamma}(t)=\gamma(t)+\lambda(t)\bv(t)$ and $\bv(t)=\pm \overline{\bw}(t)$ for all $t \in I$, 
where $\bv$ and $\bw$ are $\bt, \bn$ or $\bb$. 
We also say that $\gamma$ is a {\it $(\bv,\overline{\bw})$-Bertrand type curve} if there exists another non-degenerate regular curve $\overline{\gamma}$ such that $\gamma$ and $\overline{\gamma}$ are $(\bv,\overline{\bw})$-mates.
}
\end{definition}
We clarify the notation $\lambda \not\equiv 0$. 
Throughout this paper, $\lambda \not\equiv 0$ means that $\{t \in I | \lambda(t) \not=0\}$ is a dense subset of $I$.
Then $\lambda$ is not identically zero for any non-trivial subintervals of $I$. 
It follows that $\gamma$ and $\overline{\gamma}$ are different space curves for any non-trivial subintervals of $I$. 
Note that if $\lambda$ is constant, then $\lambda \not\equiv 0$ means that $\lambda$ is a non-zero constant.
\par
The proof of the existence conditions of $(\bt,\overline{\bn})$ and $(\bn,\overline{\bt})$-Bertrand type curves see \cite{Nakatsuyama-Takahashi}.
\begin{theorem}[\cite{Nakatsuyama-Takahashi}]\label{tn-nt-Bertrand-type}
Let $\gamma:I \to \R^3$ be non-degenerate with the arc-length parameter. 
\par
$(1)$ $\gamma$ is a $(\bt,\overline{\bn})$-Bertrand type curve if and only if  $\tau(s)=0$ and there exists a constant $c \in \R$ such that $-s+c \not=0$ for all $s \in I$. 
\par
$(2)$ $\gamma$ is an $(\bn,\overline{\bt})$-Bertrand type curve if and only if $\tau(s)=0$ and $\kappa'(s)\neq0$ for all $s \in I$.
\end{theorem}
By Theorem \ref{tn-nt-Bertrand-type}, the planar involutes and planar evolutes appear as $(\bt,\overline{\bn})$ and $(\bn,\overline{\bt})$-Bertrand type curves. 
\subsection{Framed curves}

A framed curve in the $3$-dimensional Euclidean space is a smooth space curve with a moving frame, in detail see \cite{Honda-Takahashi-2016}. 
\begin{definition}\label{framed.curve}{\rm
We say that $(\gamma,\nu_1,\nu_2):I \rightarrow \mathbb{R}^3 \times \Delta$ is a {\it framed curve} if $\dot{\gamma}(t) \cdot \nu_1(t)=0$ and $\dot{\gamma}(t) \cdot \nu_2(t)=0$ for all $t \in I$. 
We say that $\gamma:I \to \R^3$ is a {\it framed base curve} if there exists $(\nu_1,\nu_2):I \to \Delta$ such that $(\gamma,\nu_1,\nu_2)$ is a framed curve. 
}
\end{definition}

We denote $\mu(t) = \nu_1(t) \times \nu_2(t)$. 
Then $\{ \nu_1(t),\nu_2(t),\mu(t) \}$ is a moving frame along the framed base curve $\gamma(t)$ in $\R^3$ and we have the Frenet type formula,
$$
\left(
\begin{array}{c}
\dot{\nu_1}(t)\\
\dot{\nu_2}(t)\\
\dot{\mu}(t)
\end{array} \right)=
\left(
\begin{array}{ccc}
0 & \ell(t) & m(t)\\
-\ell(t) & 0 & n(t)\\
-m(t) & -n(t) & 0
\end{array}\right)
\left(
\begin{array}{c}
\nu_1(t)\\
\nu_2(t)\\
\mu(t)
\end{array}\right), 
\ \dot{\gamma}(t)=\alpha(t)\mu(t),
$$
where $\ell(t) = \dot{\nu_1}(t) \cdot \nu_2(t)$, $m(t) = \dot{\nu_1}(t) \cdot \mu(t), n(t) = \dot{\nu_2}(t) \cdot \mu(t)$ and $\alpha(t)=\dot{\gamma}(t) \cdot \mu(t)$. 
We call the mapping $(\ell,m,n,\alpha)$ {\it the curvature of the framed curve} $(\gamma,\nu_1,\nu_2)$. 
Note that $t_0$ is a singular point of $\gamma$ if and only if $\alpha(t_0) = 0$. 
\par 
The existence and uniqueness theorems for framed curves in terms of the curvatures are given by \cite{Honda-Takahashi-2016}, also see \cite{Fukunaga-Takahashi-2017}.  

\subsection{Bertrand and Mannheim curves of framed curves}

Let $(\gamma,\nu_1,\nu_2)$ and $(\overline{\gamma},\overline{\nu}_1,\overline{\nu}_2):I \to \R^3 \times \Delta$ be framed curves with the curvature $(\ell,m,n,\alpha)$ and $(\overline{\ell},\overline{m},\overline{n},\overline{\alpha})$, respectively. 
Suppose that $\gamma$ and $\overline{\gamma}$ are different curves, that is, $\gamma \not\equiv \overline{\gamma}$.

\begin{definition}\label{Bertrand-Mannheim-framed}{\rm 
We say that framed curves $(\gamma,\nu_1,\nu_2)$ and $(\overline{\gamma},\overline{\nu}_1,\overline{\nu}_2)$ are {\it Bertrand mates} (or, $(\nu_1,\overline{\nu}_1)$-mates) (respectively, {\it Mannheim mates} (or, $(\nu_1,\overline{\nu}_2)$-mates)) if there exists a smooth function $\lambda:I \to \R$ such that $\overline{\gamma}(t)=\gamma(t)+\lambda(t)\nu_1(t)$ and $\nu_1(t)=\overline{\nu}_1(t)$ (respectively, $\nu_1(t)=\overline{\nu}_2(t)$) for all $t \in I$. 
\par
We also say that $(\gamma,\nu_1,\nu_2):I \to \R^3 \times \Delta$ is a {\it Bertrand curve} (respectively, {\it Mannheim curve}) if there exists a framed curve  $(\overline{\gamma},\overline{\nu}_1,\overline{\nu}_2):I \to \R^3 \times \Delta$ such that $(\gamma,\nu_1,\nu_2)$ and $(\overline{\gamma},\overline{\nu}_1,\overline{\nu}_2)$ are Bertrand mates (respectively, Mannheim mates).
}
\end{definition}
The proof of the existence conditions of Bertrand and Mannheim curves of framed curves see \cite{Honda-Takahashi-2020}.
\begin{theorem}[\cite{Honda-Takahashi-2020}]\label{framed-Bertrand-Mannheim-equivalent}
Let $(\gamma,\nu_1,\nu_2):I \to \R^3 \times \Delta$ be a framed curve with the curvature $(\ell,m,n,\alpha)$.
\par
$(1)$ $(\gamma,\nu_1,\nu_2)$ is a Bertrand curve if and only if there exists a non-zero  constant $\lambda$ and a smooth function $\theta:I \to \R$ such that 
$
\lambda \ell(t) \cos \theta(t)-(\alpha(t)+\lambda m(t))\sin \theta(t)=0
$
for all $t \in I$.
\par
$(2)$ $(\gamma,\nu_1,\nu_2)$ is a Mannheim curve if and only if there exists a non-zero constant $\lambda$ and a smooth function $\theta:I \to \R$ such that 
$
\lambda \ell(t) \sin \theta(t)+(\alpha(t)+\lambda m(t))\cos \theta(t)=0
$
for all $t \in I$.
\end{theorem}
\subsection{$(\bv,\overline{\bw})$-Bertrand framed curves}
Let $(\gamma,\nu_1,\nu_2)$ and $(\overline{\gamma},\overline{\nu}_1,\overline{\nu}_2):I \to \R^3 \times \Delta$ be framed curves. 
\begin{definition}\label{Bertrand-type-framed}{\rm
We say that $(\gamma,\nu_1,\nu_2)$ and $(\overline{\gamma},\overline{\nu}_1,\overline{\nu}_2)$ are {\it $(\bv,\overline{\bw})$-mates} if there exists a smooth function $\lambda:I \to \R$ with $\lambda \not\equiv 0$ such that $\overline{\gamma}(t)=\gamma(t)+\lambda(t)\bv(t)$ and $\bv(t)= \overline{\bw}(t)$ for all $t \in I$, 
where $\bv$ and $\bw$ are $\nu_1, \nu_2$ or $\mu$. 
We also say that  $(\gamma,\nu_1,\nu_2)$ is a {\it $(\bv,\overline{\bw})$-Bertrand framed curve} if there exists another framed curve $(\overline{\gamma},\overline{\nu}_1,\overline{\nu}_2)$ such that $(\gamma,\nu_1,\nu_2)$ and $(\overline{\gamma},\overline{\nu}_1,\overline{\nu}_2)$ are $(\bv,\overline{\bw})$-mates.
}
\end{definition}
\par
A little bit the difference of the sign of $\overline{\bw}$ between the Definitions \ref{Bertrand-type-regular} and \ref{Bertrand-type-framed} comes from a flexibility of framed curves.
That is, if $(\gamma,\nu_1,\nu_2)$ is a framed curve, then $(\gamma,-\nu_1,\nu_2)$ and $(\gamma,\nu_1,-\nu_2)$ are also framed curves. 
Therefore, we may consider $\overline{\bw}$ up to sign.
\par
The proof of the existence conditions of $({\nu}_1, \overline{\mu})$, $({\nu}_2, \overline{\mu})$, $(\mu, \overline{\nu}_1)$ and $(\mu, \overline{\nu}_2)$-Bertrand framed curves see \cite{Nakatsuyama-Takahashi}.
\begin{theorem}[\cite{Nakatsuyama-Takahashi}]
Let $(\gamma,\nu_1,\nu_2):I \to \R^3 \times \Delta$ be a framed curve with the curvature $(\ell,m,n,\alpha)$.
\par
$(1)$ $(\gamma,\nu_1,\nu_2):I \to \R^3 \times \Delta$ is a $(\nu_1,\overline{\mu})$-Bertrand framed curve if and only if $\ell(t)=0$ and there exists a smooth function $\lambda : I\rightarrow\R$ with $\lambda \not\equiv 0$ such that $\alpha(t)+\lambda(t)m(t)=0$ for all $t\in I$. 
\par
$(2)$ $(\gamma,\nu_1,\nu_2):I \to \R^3 \times \Delta$ is a $(\nu_2,\overline{\mu})$-Bertrand framed curve if and only if $\ell(t)=0$ and there exists a smooth function $\lambda : I\rightarrow\R$ with $\lambda \not\equiv 0$ such that $\alpha(t)+\lambda(t)n(t)=0$ for all $t\in I$. 
\par
$(3)$ $(\gamma,\nu_1,\nu_2):I \to \R^3 \times \Delta$ is a $(\mu,\overline{\nu}_1)$-Bertrand framed curve if and only if there exists a smooth function $\theta:I \to \R$ such that $m(t) \cos \theta(t)-n(t)\sin \theta(t)=0$ for all $t \in I$ and $\int\alpha(t)dt\not\equiv 0$.
\par
$(4)$ $(\gamma,\nu_1,\nu_2):I \to \R^3 \times \Delta$ is a $(\mu,\overline{\nu}_2)$-Bertrand framed curve if and only if $(\gamma,\nu_1,\nu_2):I \to \R^3 \times \Delta$ is a $(\mu,\overline{\nu}_1)$-Bertrand framed curve. 
\end{theorem}
The circular evolutes (cf. \cite{Honda-Takahashi-Preprint}) (respectively, involutes) of framed curves appear as $({\nu}_1, \overline{\mu})$ and $({\nu}_2, \overline{\mu})$ (respectively, $(\mu, \overline{\nu}_1)$ and $(\mu, \overline{\nu}_2)$)-Bertrand framed curves.  

\section{$(\bv,\bw,\overline{\bx})$-Bertrand types of non-degenerate curves}

Let $\gamma$ and $\overline{\gamma}:I \to \R^3$ be non-degenerate curves with curvatures $\kappa, \overline{\kappa}$ and torsions $\tau, \overline{\tau}$. 
\begin{definition}\label{Bertrand-type-regular2}{\rm 
We say that $\gamma$ and $\overline{\gamma}$ are {\it $(\bv,\bw,\overline{\bx})$-mates} if there exists a smooth map $(\lambda,\eta):I \to \R^2$ with $(\lambda, \eta) \not\equiv (0, 0)$ such that $\overline{\gamma}(t)=\gamma(t)+\lambda(t)\bv(t)+\eta(t)\bw(t)$ and $\bv(t)\times\bw(t)=\pm\overline{\bx}(t)$ for all $t \in I$, 
where $\bv$, $\bw$ and $\bx$ are $\bt, \bn$ or $\bb$. 
We also say that $\gamma$ is a {\it $(\bv,\bw,\overline{\bx})$-Bertrand type curve} if there exists another non-degenerate regular curve $\overline{\gamma}$ such that $\gamma$ and $\overline{\gamma}$ are $(\bv,\bw,\overline{\bx})$-mates.
}
\end{definition}
Note that $(\lambda, \eta) \not\equiv (0, 0)$ means that $(\lambda, \eta)$ is not identically zero for any non-trivial subintervals of $I$. 
Then $\gamma$ and $\overline{\gamma}$ are different space curves. 

\begin{remark}\label{relation-remark}{\rm
By definition, if $\gamma$ and $\overline{\gamma}$ are $(\bv,\bw,\overline{\bx})$-mates, then $\overline{\gamma}=\gamma+\lambda\bv+\eta\bw$ and $\bv\times\bw=\pm\overline{\bx}$.
Since $\bw\times\bv=\mp\overline{\bx}$, $\gamma$ and $\overline{\gamma}=\gamma+\lambda\bw+\eta\bv$ are $(\bw,\bv,\overline{\bx})$-mates and vice versa. 
Hence, $\gamma$ is a $(\bw,\bv,\overline{\bx})$-Bertrand type curve if and only if $\gamma$ is a $(\bv,\bw,\overline{\bx})$-Bertrand type curve.
}
\end{remark}

Let $\gamma: I \to \R^3$ be a non-degenerate curve with curvature $\kappa$ and torsion $\tau$.

\begin{theorem}\label{nbb-Bertrand-type}
$(1)$ Suppose that $\tau(t)\neq0$ for all $t\in I$ and there exists a smooth function $h : I \to \R$ such that 
\begin{align}\label{evolute_non-deg}
h(t)=\frac{|\dot\gamma(t)|\tau(t)}{\kappa(t)}-\frac{d}{dt}\biggl(\frac{\dot\kappa(t)}{\kappa^2(t)}\frac{1}{|\dot\gamma(t)|\tau(t)}\biggr)\neq0.
\end{align} 
for all $t\in I$. Then $\gamma$ is an $(\bn,\bb,\overline{\bb})$-Bertrand type curve if and only if there exists a smooth map $(\lambda,\eta):I \to \R^2$ with $(\lambda, \eta) \not\equiv (0, 0)$ such that $\lambda(t)=1/\kappa(t)$ and $\eta(t)=-\dot\kappa(t)/(|\dot\gamma(t)|\kappa^2(t)\tau(t))$ for all $t \in I$. 
\par
$(2)$ Suppose that there exists a point $t\in I$ such that $\tau(t)=0$. Then $\gamma$ is not an $(\bn,\bb,\overline{\bb})$-Bertrand type curve.
\end{theorem}
\demo
$(1)$ Suppose that $\gamma$ is an $(\bn,\bb,\overline{\bb})$-Bertrand type curve. By differentiating $\overline{\gamma}(t)=\gamma(t)+\lambda(t)\bn(t)+\eta(t)\bb(t)$, we have
$
|\dot{\overline{\gamma}}(t)|\overline{\bt}(t)=|\dot{\gamma}(t)|(1-\lambda(t)\kappa(t))\bt (t)+(\dot\lambda(t)-|\dot{\gamma}(t)|\eta(t)\tau(t))\bn(t)+(|\dot{\gamma}(t)|\lambda(t)\tau(t)+\dot\eta(t))\bb(t).
$
Since $\bn(t)\times\bb(t)=\bt(t)=\pm\overline{\bb}(t)$, we have $1-\lambda(t)\kappa(t)=0$ for all $t\in I$. Therefore, we have $\lambda(t)=1/\kappa(t)$.
If $\bt(t)=\overline{\bb}(t)$, there exists a smooth function $\theta:I \to \R$ such that
$$
\left(
\begin{array}{c}
\overline{\bt}(t)\\
\overline{\bn}(t)
\end{array}
\right)
=
\left(
\begin{array}{cc}
\cos\theta(t) & -\sin\theta(t)\\
\sin\theta(t) & \cos\theta(t)
\end{array}
\right)
\left(
\begin{array}{c}
\bn(t)\\
\bb(t)
\end{array}
\right).
$$
Then, we have $|\dot{\overline{\gamma}}(t)|\cos\theta(t)=\dot\lambda(t)-|\dot{\gamma}(t)|\eta(t)\tau(t)$ and $|\dot{\overline{\gamma}}(t)|\sin\theta(t)=-(|\dot{\gamma}(t)|\lambda(t)\tau(t)+\dot\eta(t))$.
By differentiating $\overline{\bt}(t)=\cos\theta(t)\bn(t)-\sin\theta(t)\bb(t)$, then $\cos\theta(t)=0$ and $\sin\theta(t)=\pm1$. 
Thus, we have $\dot\lambda(t)-|\dot{\gamma}(t)|\eta(t)\tau(t)=0$ for all $t\in I$. 
Therefore, we have 
$$
\eta(t)=\frac{\dot{\lambda}(t)}{|\dot{\gamma}(t)|\tau(t)}=-\frac{\dot{\kappa}(t)}{\vert \dot{\gamma}(t) \vert \kappa^2(t) \tau(t)}
$$ 
and $|\dot{\overline{\gamma}}(t)|=\mp(|\dot{\gamma}(t)|\lambda(t)\tau(t)+\dot\eta(t))=\mp(|\dot{{\gamma}}(t)|{\tau}(t)/\kappa(t)-(d/dt)({\dot{\kappa}(t)}/({\vert \dot{\gamma}(t) \vert \kappa^2(t) \tau(t)})))$.
Therefore, we have $|\dot{\overline{\gamma}}(t)|=\mp h(t)\neq0$. 
\par
If $\bt(t)=-\overline{\bb}(t)$, we also obtain the same result as $\bt(t)=\overline{\bb}(t)$.
\par
Conversely, suppose that $\lambda(t)=1/\kappa(t)$ and $\eta(t)=-{\dot{\kappa}(t)}/({\vert \dot{\gamma}(t) \vert \kappa^2(t) \tau(t)})$ for all $t \in I$. 
By a direct calculation, we have
\begin{align*}
\dot{\overline{\gamma}}(t)&=h(t)\bb(t), \ \ddot{\overline{\gamma}}(t)=\dot{h}(t)\bb(t)-h(t)|\dot{\gamma}(t)|\tau(t)\bn(t), \ \dot{\overline{\gamma}}(t)\times\ddot{\overline{\gamma}}(t)=h^2(t)|\dot{\gamma}(t)|\tau(t)\bt(t)\neq0.
\end{align*}
It follows that $\overline{\gamma}$ is a non-degenerate curve. Moreover, we have 
\begin{align*}
{\overline{\bb}}(t)=\frac{\dot{\overline{\gamma}}(t)\times\ddot{\overline{\gamma}}(t)}{|\dot{\overline{\gamma}}(t)\times\ddot{\overline{\gamma}}(t)|}=\frac{h^2(t)|\dot{\gamma}(t)|\tau(t)}{h^2(t)|\dot{\gamma}(t)||\tau(t)|}\bt(t)=\pm\bt(t).
\end{align*}
It follows that $\gamma$ is an $(\bn,\bb,\overline{\bb})$-Bertrand type curve.
\par
$(2)$ Suppose that there exists a point $t\in I$ such that $\tau(t)=0$ and $\gamma$ is an $(\bn,\bb,\overline{\bb})$-Bertrand type curve.
By the same method of $(1)$, we have $\dot\gamma(t)\times\ddot\gamma(t)=h^2(t)|\dot\gamma(t)|\tau(t)\bt(t)$.
Then, we have $\dot\gamma(t)\times\ddot\gamma(t)=0$ at the point. 
Therefore, $\overline{\gamma}$ is not non-degenerate and $\gamma$ is not an $(\bn,\bb,\overline{\bb})$-Bertrand type curve.
\enD
\begin{remark}\label{nbb-Bertrand-type_curvature}{\rm
If $\gamma$ is an $(\bn,\bb,\overline{\bb})$-Bertrand type curve, then $\overline{\gamma}$ is an evolute of $\gamma$  (cf. \cite{Blaschke, Bruce-Giblin, Gibson,Gray}). 
If $h>0$ (respectively, $h<0$), the moving frame of $\overline{\gamma}$ is given by $\{\overline{\bt},\overline{\bn},\overline{\bb}\}=\{\bb,\mp\bn,\pm\bt\}$ (respectively, $\{\overline{\bt},\overline{\bn},\overline{\bb}\}=\{-\bb,\pm\bn,\pm\bt\}$). 
Moreover, the curvature $\overline{\kappa}$ and the torsion $\overline{\tau}$ of $\overline{\gamma}$ are given by
$
\overline{\kappa}(t)={|\dot\gamma(t)||\tau(t)|}/{|h(t)|}, \overline{\tau}(t)={|\dot\gamma(t)|\kappa(t)}/{h(t)}.
$
}
\end{remark}

\begin{theorem}\label{tnt-Bertrand-type}
$(1)$ Suppose that $\tau(t)\neq0$ for all $t\in I$.
Then $\gamma$ is a $(\bt,\bn,\overline{\bt})$-Bertrand type curve if and only if 
there exists a smooth map $(\lambda,\eta):I \to \R^2$ with $\eta(t) \neq 0$ such that
\begin{align}\label{tnt-Bertrand-condition}
\begin{cases}
|\dot\gamma(t)|+\dot\lambda(t)-\eta(t)|\dot\gamma(t)|\kappa(t)=0,\\
\lambda(t)|\dot\gamma(t)|\kappa(t)+\dot\eta(t)=0
\end{cases}
\end{align}
for all $t \in I$.
\par 
$(2)$ Suppose that there exists a point $t\in I$ such that $\tau(t)=0$. Then $\gamma$ is not a $(\bt,\bn,\overline{\bt})$-Bertrand type curve.
\end{theorem}
\demo
$(1)$ Suppose that $\gamma$ is a $(\bt,\bn,\overline{\bt})$-Bertrand type curve. By differentiating $\overline{\gamma}(t)=\gamma(t)+\lambda(t)\bt(t)+\eta(t)\bn(t)$, we have
$
|\dot{\overline{\gamma}}(t)|\overline{\bt}(t)=(|\dot\gamma(t)|+\dot\lambda(t)-\eta(t)|\dot\gamma(t)|\kappa(t))\bt(t)+(\lambda(t)|\dot\gamma(t)|\kappa(t)+\dot\eta(t))\bn(t)+\eta(t)|\dot\gamma(t)|\tau(t)\bb(t).
$
Since $\bt(t)\times\bn(t)=\bb(t)=\pm\overline{\bt}(t)$, we have 
$$
|\dot{\overline{\gamma}}(t)|=\pm\eta(t)|\dot\gamma(t)|\tau(t), \ |\dot\gamma(t)|+\dot\lambda(t)-\eta(t)|\dot\gamma(t)|\kappa(t)=0, \
\lambda(t)|\dot\gamma(t)|\kappa(t)+\dot\eta(t)=0
$$ 
for all $t \in I$. Therefore, we have $\eta(t)\neq0$ for all $t \in I$.
\par
Conversely, suppose that there exists a smooth map $(\lambda,\eta):I \to \R^2$ with $\eta(t) \neq 0$ such that condition \eqref{tnt-Bertrand-condition} satisfies for all $t \in I$.
By a direct calculation, we have
\begin{align*}
\dot{\overline{\gamma}}(t)&=\eta(t)|\dot\gamma(t)|\tau(t)\bb(t), \ 
\ddot{\overline{\gamma}}(t)=\frac{d}{dt}(\eta(t)|\dot\gamma(t)|\tau(t))\bb(t)-\eta(t)|\dot\gamma(t)|^2\tau^2(t)\bn(t), \\ 
\dot{\overline{\gamma}}(t)&\times\ddot{\overline{\gamma}}(t)=\eta^2(t)|\dot\gamma(t)|^3\tau^3(t)\bt(t)\neq0.
\end{align*}
It follows that $\overline{\gamma}$ is a non-degenerate curve. Moreover, we have 
\begin{align*}
{\overline{\bt}}(t)=\frac{\dot{\overline{\gamma}}(t)}{|\dot{\overline{\gamma}}(t)|}=\frac{\eta(t)|\dot\gamma(t)|\tau(t)}{|\eta(t)||\dot\gamma(t)||\tau(t)|}\bb(t)=\pm\bb(t).
\end{align*}
It follows that $\gamma$ is a $(\bt,\bn,\overline{\bt})$-Bertrand type curve.
\par
$(2)$ Suppose that there exists a point $t\in I$ such that $\tau(t)=0$ and $\gamma$ is a $(\bt,\bn,\overline{\bt})$-Bertrand type curve.
By the same method of $(1)$, we have $\dot\gamma(t)\times\ddot\gamma(t)=\eta^2(t)|\dot\gamma(t)|^3\tau^3(t)\bt(t)$.
Then, we have $\dot\gamma(t)\times\ddot\gamma(t)=0$ at the point. 
Therefore, $\overline{\gamma}$ is not non-degenerate and $\gamma$ is not a $(\bt,\bn,\overline{\bt})$-Bertrand type curve.
\enD

\begin{remark}\label{tnt-Bertrand-type_curvature}{\rm
If $\gamma$ is a $(\bt,\bn,\overline{\bt})$-Bertrand type curve, then we may consider $\overline{\gamma}$ is one of involutes of $\gamma$. 
If $\tau>0$ (respectively, $\tau<0$), the moving frame of $\overline{\gamma}$ is given by $\{\overline{\bt},\overline{\bn},\overline{\bb}\}=\{\pm\bb,\mp\bn,\bt\}$ (respectively, $\{\overline{\bt},\overline{\bn},\overline{\bb}\}=\{\pm\bb,\pm\bn,-\bt\}$). 
Moreover, the curvature $\overline{\kappa}$ and the torsion $\overline{\tau}$ of $\overline{\gamma}$ are given by $\overline{\kappa}(t)={1}/{|\eta(t)|},  \overline{\tau}(t)={\kappa(t)}/({\eta(t)\tau(t)})$. 
}
\end{remark}

\begin{theorem}\label{nbt-Bertrand-type}
$\gamma$ is an $(\bn,\bb,\overline{\bt})$-Bertrand type curve if and only if there there exists a smooth map $(\lambda,\eta):I \to \R^2$ with $(\lambda,\eta) \not\equiv (0,0)$ such that
\begin{align}\label{nbt-Bertrand-type-condition}
\begin{cases}
1-\lambda(t)\kappa(t)\neq0,\ \dot\lambda(t)-\eta(t)|\dot\gamma(t)|\tau(t)=0,\\
\lambda(t)|\dot\gamma(t)|\tau(t)+\dot\eta(t)=0
\end{cases}
\end{align}
for all $t \in I$.
\end{theorem}
\demo
Suppose that $\gamma$ is an $(\bn,\bb,\overline{\bt})$-Bertrand type curve. By differentiating $\overline{\gamma}(t)=\gamma(t)+\lambda(t)\bn(t)+\eta(t)\bb(t)$, we have
$
|\dot{\overline{\gamma}}(t)|\overline{\bt}(t)=|\dot\gamma(t)|(1-\lambda(t)\kappa(t))\bt(t)+(\dot\lambda(t)-\eta(t)|\dot\gamma(t)|\tau(t))\bn(t)+(\lambda(t)|\dot\gamma(t)|\tau(t)+\dot\eta(t))\bb(t).
$
Since $\bn(t)\times\bb(t)=\bt(t)=\pm\overline{\bt}(t)$, we have 
$$
|\dot{\overline{\gamma}}(t)|=\pm|\dot\gamma(t)|(1-\lambda(t)\kappa(t)), \
\dot\lambda(t)-\eta(t)|\dot\gamma(t)|\tau(t)=0, \ \lambda(t)|\dot\gamma(t)|\tau(t)+\dot\eta(t)=0
$$ 
for all $t \in I$. 
Therefore, we have $1-\lambda(t)\kappa(t)\neq0$.
\par
Conversely, suppose that there exists a smooth map $(\lambda,\eta):I \to \R^2$ with $(\lambda, \eta) \not\equiv (0, 0)$ such that condition \eqref{nbt-Bertrand-type-condition} satisfies for all $t \in I$.
By a direct calculation, we have
\begin{align*}
\dot{\overline{\gamma}}(t)&=|\dot\gamma(t)|(1-\lambda(t)\kappa(t))\bt(t), \\ 
\ddot{\overline{\gamma}}(t)&=\frac{d}{dt}\left(|\dot\gamma(t)|(1-\lambda(t)\kappa(t))\right)\bt(t)+|\dot\gamma(t)|^2\kappa(t)(1-\lambda(t)\kappa(t))\bn(t), \\ 
\dot{\overline{\gamma}}(t)&\times\ddot{\overline{\gamma}}(t)=|\dot\gamma(t)|^3\kappa(t)(1-\lambda(t)\kappa(t))^2\bb(t)\neq0.
\end{align*}
It follows that $\overline{\gamma}$ is a non-degenerate curve. Moreover, we have 
\begin{align*}
{\overline{\bt}}(t)=\frac{\dot{\overline{\gamma}}(t)}{|\dot{\overline{\gamma}}(t)|}=\frac{|\dot\gamma(t)|(1-\lambda(t)\kappa(t))}{|\dot\gamma(t)||1-\lambda(t)\kappa(t)|}\bt(t)=\pm\bt(t).
\end{align*}
It follows that $\gamma$ is an $(\bn,\bb,\overline{\bt})$-Bertrand type curve.
\enD
\begin{remark}\label{nbt-Bertrand-type_curvature}{\rm
The moving frame of $\overline{\gamma}$ is given by $\{\overline{\bt},\overline{\bn},\overline{\bb}\}=\{\pm\bt,\bn,\pm\bb\}$. 
The curvature $\overline{\kappa}$ and the torsion $\overline{\tau}$ of $\overline{\gamma}$ are given by $\overline{\kappa}(t)={\kappa(t)}/{(1-\lambda(t)\kappa(t))}$ and $\overline{\tau}(t)={\tau(t)}/{(1-\lambda(t)\kappa(t))}$. 
}
\end{remark}

\section{$(\bv,\bw,\overline{\bx})$-Bertrand framed curves}

Let $(\gamma,\nu_1,\nu_2)$ and $(\overline{\gamma},\overline{\nu}_1,\overline{\nu}_2):I \to \R^3 \times \Delta$ be framed curves with curvatures $(\ell,m,n,\alpha)$ and $(\overline{\ell},\overline{m},\overline{n},\overline{\alpha})$. 
\begin{definition}\label{Bertrand-type-framed2}{\rm
We say that $(\gamma,\nu_1,\nu_2)$ and $(\overline{\gamma},\overline{\nu}_1,\overline{\nu}_2)$ are {\it $(\bv,\bw,\overline{\bx})$-mates} if there exists a smooth map $(\lambda,\eta):I \to \R^2$ with $(\lambda,\eta) \not\equiv (0,0)$ such that $\overline{\gamma}(t)=\gamma(t)+\lambda(t)\bv(t)+\eta(t)\bw(t)$ and $\bv(t)\times\bw(t)= \overline{\bx}(t)$ for all $t \in I$, 
where $\bv$, $\bw$ and $\bx$ are $\nu_1, \nu_2$ or $\mu$. 
We also say that $(\gamma,\nu_1,\nu_2)$ is a {\it $(\bv,\bw,\overline{\bx})$-Bertrand framed curve} if there exists a framed curve $(\overline{\gamma},\overline{\nu}_1,\overline{\nu}_2)$ such that $(\gamma,\nu_1,\nu_2)$ and $(\overline{\gamma},\overline{\nu}_1,\overline{\nu}_2)$ are $(\bv,\bw,\overline{\bx})$-mates.
}
\end{definition}
Note that $(\lambda, \eta) \not\equiv (0, 0)$ means that $(\lambda, \eta)$ is not identically zero for any non-trivial subintervals of $I$. 
Then $\gamma$ and $\overline{\gamma}$ are different framed base curves.
A little bit the difference of the sign of $\overline{\bx}$ between the Definitions \ref{Bertrand-type-regular2} and \ref{Bertrand-type-framed2} comes from a flexibility of framed curves.
That is, if $(\gamma,\nu_1,\nu_2)$ is a framed curve, then $(\gamma,\nu_2,\nu_1)$, $(\gamma,-\nu_1,\nu_2)$ and $(\gamma,\nu_1,-\nu_2)$ are also framed curves. 
\par
By definition, we have the following.
\begin{proposition}\label{frame-change}{\rm
Suppose that $(\gamma,\nu_1,\nu_2):I \to \R^3 \times \Delta$ is a $(\bv,\bw,\overline{\bx})$-Bertrand framed curve. 
\par
$(1)$ $(\gamma,\nu_1,\nu_2)$ is also a $(\bw,\bv,-\overline{\bx})$-Bertrand framed curve.
\par
$(2)$ $(\gamma,\nu_1,\nu_2)$ is also a $(-\bv,\bw,-\overline{\bx})$-Bertrand framed curve.
\par
$(3)$ $(\gamma,\nu_1,\nu_2)$ is also a $(\bv,-\bw,-\overline{\bx})$-Bertrand framed curve.
}
\end{proposition}
Let $(\gamma,\nu_1,\nu_2):I \to \R^3 \times \Delta$ be a framed curve with curvature $(\ell,m,n,\alpha)$. 
\begin{theorem}\label{n1n2n2-Bertrand-type}
$(\gamma,\nu_1,\nu_2):I \to \R^3 \times \Delta$ is a $(\nu_1,\nu_2,\overline{\nu}_2)$-Bertrand framed curve if and only if there exist a smooth map $(\lambda,\eta):I \to \R^2$ with $(\lambda, \eta) \not\equiv (0, 0)$ and a smooth function $\theta:I \to \R$ such that
\begin{align}\label{n1n2n2-Bertrand-type-condition}
\begin{cases}
\alpha(t)+\lambda(t) m(t)+\eta(t) n(t)=0,\\
(\dot\lambda(t)-\eta(t)\ell(t))\sin\theta(t)+(\lambda(t)\ell(t)+\dot\eta(t))\cos\theta(t)=0
\end{cases}
\end{align}
for all $t \in I$.
\end{theorem}
\demo
Suppose that $(\gamma,\nu_1,\nu_2):I \to \R^3 \times \Delta$ is a $(\nu_1,\nu_2,\overline{\nu}_2)$-Bertrand framed curve. 
By differentiating $\overline{\gamma}(t)=\gamma(t)+\lambda(t)\nu_1(t)+\eta(t)\nu_2(t)$, we have 
$
\overline{\alpha}(t) \overline{\mu}(t)=(\dot\lambda(t)-\eta(t)\ell(t))\nu_1(t)+(\lambda(t)\ell(t)+\dot\eta(t))\nu_2(t)+(\alpha(t)+\lambda(t)m(t)+\eta(t)n(t))\mu(t)
$ for all $t \in I$. 
Since $\nu_1(t)\times\nu_2(t)=\mu(t)=\overline{\nu}_2(t)$, we have $\alpha(t)+\lambda(t)m(t)+\eta(t)n(t)=0$ for all $t \in I$.
Moreover, there exists a smooth function $\theta:I \to \R$ such that 
$$
\begin{pmatrix}
\overline{\mu}(t) \\
\overline{\nu}_1(t)
\end{pmatrix}
=
\begin{pmatrix}
\cos \theta(t) & -\sin \theta(t) \\
\sin \theta(t) & \cos \theta(t)
\end{pmatrix}
\begin{pmatrix}
{\nu}_1(t) \\
{\nu}_2(t)
\end{pmatrix}.
$$
Then, we have $\overline{\alpha}(t)\sin \theta(t)=-(\lambda(t)\ell(t)+\dot\eta(t))$ and $\overline{\alpha}(t)\cos(t)=\dot\lambda(t)-\eta(t)\ell(t)$.
It follows that $(\dot\lambda(t)-\eta(t)\ell(t))\sin\theta(t)+(\lambda(t)\ell(t)+\dot\eta(t))\cos\theta(t)=0$ for all $t \in I$. 
\par
Conversely, suppose that there exist a smooth map $(\lambda,\eta):I \to \R^2$ with $(\lambda, \eta) \not\equiv (0, 0)$ and a smooth function $\theta:I \to \R$ such that condition (\ref{n1n2n2-Bertrand-type-condition}) satisfies. 
Let $(\overline{\gamma},\overline{\nu}_1,\overline{\nu}_2):I \to \R^3 \times \Delta$ be 
$
\overline{\gamma}(t)=\gamma(t)+\lambda(t)\nu_1(t)+\eta(t)\nu_2(t), \ \overline{\nu}_1(t)=\sin\theta(t)\nu_1(t)+\cos\theta(t)\nu_2(t)$
 and $
\overline{\nu}_2(t)=\mu(t).
$
Since $\dot{\overline{\gamma}}(t)=(\dot\lambda(t)-\eta(t)\ell(t))\nu_1(t)+(\lambda(t)\ell(t)+\dot\eta(t))\nu_2(t)$, then we have
$\dot{\overline{\gamma}}(t)\cdot\overline{\nu}_1(t)=\dot{\overline{\gamma}}(t)\cdot\overline{\nu}_2(t)=0$. 
It follows that $(\overline{\gamma},\overline{\nu}_1,\overline{\nu}_2)$ is a framed curve.
Moreover, we have $\overline{\nu}_2(t)=\mu(t)=\nu_1(t)\times\nu_2(t)$.
Therefore, $(\gamma,\nu_1,\nu_2)$ is a $(\nu_1,\nu_2,\overline{\nu}_2)$-Bertrand framed curve.
\enD
\begin{remark}\label{n1n2n2-Bertrand-type_evolute}{\rm
If $(\gamma,\nu_1,\nu_2)$ is a $(\nu_1,\nu_2,\overline{\nu}_2)$-Bertrand framed curve, then $\overline{\gamma}$ is an evolute of $(\gamma,\nu_1,\nu_2)$ (cf. \cite{Honda-Takahashi-2020-evolute}). 
Moreover, $(\gamma,\nu_1,\nu_2)$ is a $(\nu_1,\nu_2,\overline{\nu}_1)$-Bertrand framed curve if and only if $(\gamma,\nu_1,\nu_2)$ is a $(\nu_1,\nu_2,\overline{\nu}_2)$-Bertrand framed curve.
}
\end{remark}
\begin{proposition}\label{n1n2n2-Bertrand-type_curvature}
Suppose that $(\gamma,\nu_1,\nu_2)$ and $(\overline{\gamma},\overline{\nu}_1,\overline{\nu}_2):  I \to \R^3 \times \Delta$ are {$(\nu_1,\nu_2,\overline{\nu}_2)$-mates}, where $\overline{\gamma}(t)=\gamma(t)+\lambda(t)\nu_1(t)+\eta(t)\nu_2(t), \overline{\nu}_1(t)=\sin\theta(t)\nu_1(t)+\cos\theta(t)\nu_2(t)$ and $\overline{\nu}_2(t)=\mu(t)$.  
Then the curvature $(\overline{\ell}, \overline{m}, \overline{n}, \overline{\alpha})$ of $(\overline{\gamma},\overline{\nu}_1,\overline{\nu}_2)$ is given by
\begin{align*}
\overline{\ell}(t)&=m(t)\sin\theta(t)+n(t)\cos\theta(t), \ \overline{m}(t)=\dot\theta(t)-\ell(t), \\
\overline{n}(t)&=n(t)\sin\theta(t)-m(t)\cos\theta(t), \ \overline{\alpha}(t)=(\dot\lambda(t)-\eta(t)\ell(t))\cos\theta(t)-(\lambda(t)\ell(t)+\dot\eta(t))\sin\theta(t). 
\end{align*}
\end{proposition}
\demo
By differentiating $\overline{\nu}_1(t)=\sin\theta(t)\nu_1(t)+\cos\theta(t)\nu_2(t)$ and $\overline{\mu}(t)=\cos\theta(t)\nu_1(t)-\sin\theta(t)\nu_2(t)$, we have 
\begin{align*}
\overline{\ell}(t)\overline{\nu}_2(t)+\overline{m}(t)\overline{\mu}(t)&=(\dot\theta(t)-\ell(t))\overline{\mu}(t)+(m(t)\sin\theta(t)+n(t)\cos\theta(t))\overline{\nu}_2(t),\\
-\overline{m}(t)\overline{\nu}_1(t)-\overline{n}(t)\overline{\nu}_2(t)&=-(\dot\theta(t)-\ell(t))\overline{\nu}_1(t)+(m(t)\cos\theta(t)-n(t)\sin\theta(t))\overline{\nu}_2(t).
\end{align*}
Then, we have 
$$
\overline{\ell}(t)=m(t)\sin\theta(t)+n(t)\cos\theta(t), \ 
\overline{m}(t)=\dot\theta(t)-\ell(t), \ 
\overline{n}(t)=n(t)\sin\theta(t)-m(t)\cos\theta(t).
$$
By Theorem \ref{n1n2n2-Bertrand-type}, we have $\overline{\alpha}(t)=(\dot\lambda(t)-\eta(t)\ell(t))\cos\theta(t)-(\lambda(t)\ell(t)+\dot\eta(t))\sin\theta(t)$.
\enD

\begin{theorem}\label{mun1mu-Bertrand-type}
$(\gamma,\nu_1,\nu_2):I \to \R^3 \times \Delta$ is a $(\mu,\nu_1,\overline{\mu})$-Bertrand framed curve if and only if there exists a smooth map $(\lambda,\eta):I \to \R^2$ with $(\lambda, \eta) \not\equiv (0, 0)$ such that
\begin{align}\label{mun1mu-Bertrand-type-condition}
\begin{cases}
-\lambda(t)m(t)+\dot\eta(t)=0, \\
\alpha(t)+\dot\lambda(t)+\eta(t)m(t)=0 
\end{cases}
\end{align}
for all $t \in I$.
\end{theorem}
\demo
Suppose that $(\gamma,\nu_1,\nu_2):I \to \R^3 \times \Delta$ is a $(\mu,\nu_1,\overline{\mu})$-Bertrand framed curve. 
By differentiating $\overline{\gamma}(t)=\gamma(t)+\lambda(t)\mu(t)+\eta(t)\nu_1(t)$, we have 
$
\overline{\alpha}(t) \overline{\mu}(t)=(-\lambda(t)m(t)+\dot\eta(t))\nu_1(t)+(-\lambda(t)n(t)+\eta(t)\ell(t))\nu_2(t)+(\alpha(t)+\dot\lambda(t)+\eta(t)m(t))\mu(t)
$ for all $t \in I$. 
Since $\mu(t)\times\nu_1(t)=\nu_2(t)=\overline{\mu}(t)$, we have 
$$
\overline{\alpha}(t)=-\lambda(t)n(t)+\eta(t)\ell(t), \
-\lambda(t)m(t)+\dot\eta(t)=0, \
\alpha(t)+\dot\lambda(t)+\eta(t)m(t)=0
$$ 
for all $t \in I$. 
\par
Conversely, suppose that there exists a smooth map $(\lambda,\eta):I \to \R^2$ with $(\lambda, \eta) \not\equiv (0, 0)$ such that condition (\ref{mun1mu-Bertrand-type-condition}) satisfies. 
Let $(\overline{\gamma},\overline{\nu}_1,\overline{\nu}_2):I \to \R^3 \times \Delta$ be 
$
\overline{\gamma}(t)=\gamma(t)+\lambda(t)\mu(t)+\eta(t)\nu_1(t), \ \overline{\nu}_1(t)=\cos\theta(t)\mu(t)-\sin\theta(t)\nu_1(t)
$ and 
$\overline{\nu}_2(t)=\sin\theta(t)\mu(t)+\cos\theta(t)\nu_1(t),
$
where $\theta:I\to\R$ is a smooth function.
Since $\dot{\overline{\gamma}}(t)=(-\lambda(t)n(t)+\eta(t)\ell(t))\nu_2(t)$, then we have
$\dot{\overline{\gamma}}(t)\cdot\overline{\nu}_1(t)=\dot{\overline{\gamma}}(t)\cdot\overline{\nu}_2(t)=0$. 
It follows that $(\overline{\gamma},\overline{\nu}_1,\overline{\nu}_2)$ is a framed curve.
Moreover, by a direct calculation, we have $\overline{\mu}(t)=\overline{\nu}_1(t)\times\overline{\nu}_2(t)=(\cos^2\theta(t)+\sin^2\theta(t))\nu_2(t)=\nu_2(t)=\mu(t)\times\nu_1(t)$.
Therefore, $(\gamma,\nu_1,\nu_2)$ is a $(\mu,\nu_1,\overline{\mu})$-Bertrand framed curve.
\enD
\begin{remark}\label{mun1mu-Bertrand-type_involute}{\rm
If $(\gamma,\nu_1,\nu_2)$ is a $(\mu,\nu_1,\overline{\mu})$-Bertrand framed curve, then we may consider $\overline{\gamma}$ is one of involute of $(\gamma,\nu_1,\nu_2)$. 
Moreover, $(\gamma,\nu_1,\nu_2)$ is a $(\mu,\nu_2,\overline{\mu})$-Bertrand framed curve if and only if $(\gamma,\nu_1,\nu_2)$ is a $(\mu,\nu_1,\overline{\mu})$-Bertrand framed curve.
}
\end{remark}
By the same method of Proposition \ref{n1n2n2-Bertrand-type_curvature}, we have the following.
\begin{proposition}\label{mun1mu-Bertrand-type_curvature}
Suppose that $(\gamma,\nu_1,\nu_2)$ and $(\overline{\gamma},\overline{\nu}_1,\overline{\nu}_2):I \to \R^3 \times \Delta$ are {$(\mu,\nu_1,\overline{\mu})$-mates}, where $\overline{\gamma}(t)=\gamma(t)+\lambda(t)\mu(t)+\eta(t)\nu_1(t), \overline{\nu}_1(t)=\cos\theta(t)\mu(t)-\sin\theta(t)\nu_1(t)$ and $\overline{\nu}_2(t)=\sin\theta(t)\mu(t)+\cos\theta(t)\nu_1(t)$. 
Then the curvature $(\overline{\ell}, \overline{m}, \overline{n}, \overline{\alpha})$ of $(\overline{\gamma},\overline{\nu}_1,\overline{\nu}_2)$ is given by
\begin{align*}
  \overline{\ell}(t)&=-\dot\theta(t)-m(t), \ \overline{m}(t)=-(\ell(t)\sin\theta(t)+n(t)\cos\theta(t)), \\
  \overline{n}(t)&=\ell(t)\cos\theta(t)-n(t)\sin\theta(t), \ \overline{\alpha}(t)=-\lambda(t)n(t)+\eta(t)\ell(t). 
\end{align*}
\end{proposition}

\begin{theorem}\label{mun1n2-Bertrand-type}
$(\gamma,\nu_1,\nu_2):I \to \R^3 \times \Delta$ is a $(\mu,\nu_1,\overline{\nu}_2)$-Bertrand framed curve if and only if there exist a smooth map $(\lambda,\eta):I \to \R^2$ with $(\lambda, \eta) \not\equiv (0, 0)$ and a smooth function $\theta:I\to\R$ such that
\begin{align}\label{mun1n2-Bertrand-type-condition}
\begin{cases}
-\lambda(t)n(t)+\eta(t)\ell(t)=0, \\ 
(\alpha(t)+\dot\lambda(t)+\eta(t) m(t))\sin\theta(t)+(-\lambda(t) m(t)+\dot\eta(t))\cos\theta(t)=0
\end{cases}
\end{align}
for all $t \in I$.
\end{theorem}
\demo
Suppose that $(\gamma,\nu_1,\nu_2):I \to \R^3 \times \Delta$ is a $(\mu,\nu_1,\overline{\nu}_2)$-Bertrand framed curve. 
By differentiating $\overline{\gamma}(t)=\gamma(t)+\lambda(t)\mu(t)+\eta(t)\nu_1(t)$, we have 
$
\overline{\alpha}(t) \overline{\mu}(t)=(-\lambda(t)m(t)+\dot\eta(t))\nu_1(t)+(-\lambda(t)n(t)+\eta(t)\ell(t))\nu_2(t)+(\alpha(t)+\dot\lambda(t)+\eta(t)m(t))\mu(t)
$ for all $t \in I$. 
Since $\mu(t)\times\nu_1(t)=\nu_2(t)=\overline{\nu}_2(t)$, we have $-\lambda(t)n(t)+\eta(t)\ell(t)=0$ for all $t \in I$. 
Moreover, there exists a smooth function $\theta:I \to \R$ such that 
$$
\begin{pmatrix}
\overline{\mu}(t) \\
\overline{\nu}_1(t)
\end{pmatrix}
=
\begin{pmatrix}
\cos \theta(t) & -\sin \theta(t) \\
\sin \theta(t) & \cos \theta(t)
\end{pmatrix}
\begin{pmatrix}
{\mu}(t) \\
{\nu}_1(t)
\end{pmatrix}.
$$
Then, we have $\overline{\alpha}(t)\sin \theta(t)=\lambda(t)m(t)-\dot\eta(t)$ and $\overline{\alpha}(t)\cos\theta(t)=\alpha(t)+\dot\lambda(t)+\eta(t)m(t)$.
It follows that $(\alpha(t)+\dot\lambda(t)+\eta(t)m(t))\sin\theta(t)+(-\lambda(t) m(t)+\dot\eta(t))\cos\theta(t)=0$ for all $t \in I$. 
\par
Conversely, suppose that there exist a smooth map $(\lambda,\eta):I \to \R^2$ with $(\lambda, \eta) \not\equiv (0, 0)$ and a smooth function $\theta:I \to \R$ such that condition (\ref{mun1n2-Bertrand-type-condition}) satisfies. 
Let $(\overline{\gamma},\overline{\nu}_1,\overline{\nu}_2):I \to \R^3 \times \Delta$ be 
$
\overline{\gamma}(t)=\gamma(t)+\lambda(t)\mu(t)+\eta(t)\nu_1(t), \ \overline{\nu}_1(t)=\sin\theta(t)\mu(t)+\cos\theta(t)\nu_1(t)
$ and 
$\overline{\nu}_2(t)=\nu_2(t),
$
where $\theta:I\to\R$ is a smooth function. 
Since $\dot{\overline{\gamma}}(t)=(\lambda(t)m(t)-\dot\eta(t))\nu_1(t)+(\alpha(t)+\dot\lambda(t)+\eta(t)m(t))\mu(t)$, we have
$\dot{\overline{\gamma}}(t)\cdot\overline{\nu}_1(t)=\dot{\overline{\gamma}}(t)\cdot\overline{\nu}_2(t)=0$. 
It follows that $(\overline{\gamma},\overline{\nu}_1,\overline{\nu}_2)$ is a framed curve.
Moreover, we have $\overline{\nu}_2(t)=\nu_2(t)=\mu(t)\times\nu_1(t)$.
Therefore, $(\gamma,\nu_1,\nu_2)$ is a $(\mu,\nu_1,\overline{\nu}_2)$-Bertrand framed curve.
\enD

\begin{remark}\label{mun1n2-Bertrand-type-equivalent}{\rm
$(\gamma,\nu_1,\nu_2)$ is a $(\mu,\nu_2,\overline{\nu}_2)$-Bertrand framed curve if and only if $(\gamma,\nu_1,\nu_2)$ is a $(\mu,\nu_1,\overline{\nu}_2)$-Bertrand framed curve.
Moreover, $(\gamma,\nu_1,\nu_2)$ is a $(\mu,\nu_1,\overline{\nu}_1)$-Bertrand framed curve if and only if $(\gamma,\nu_1,\nu_2)$ is a $(\mu,\nu_1,\overline{\nu}_2)$-Bertrand framed curve.
Hence, $(\gamma,\nu_1,\nu_2)$ is a $(\mu,\nu_2,\overline{\nu}_1)$-Bertrand framed curve if and only if $(\gamma,\nu_1,\nu_2)$ is a $(\mu,\nu_1,\overline{\nu}_2)$-Bertrand framed curve.
}
\end{remark}
By the same method of Proposition \ref{n1n2n2-Bertrand-type_curvature}, we have the following.
\begin{proposition}\label{mun1n2-Bertrand-type_curvature}
Suppose that $(\gamma,\nu_1,\nu_2)$ and $(\overline{\gamma},\overline{\nu}_1,\overline{\nu}_2):  I \to \R^3 \times \Delta$ are {$(\mu,\nu_1,\overline{\nu}_2)$-mates}, where $\overline{\gamma}(t)=\gamma(t)+\lambda(t)\mu(t)+\eta(t)\nu_1(t), \overline{\nu}_1(t)=\sin\theta(t)\mu(t)+\cos\theta(t)\nu_1(t)$ and $\overline{\nu}_2(t)=\nu_2(t)$. 
Then the curvature $(\overline{\ell}, \overline{m}, \overline{n}, \overline{\alpha})$ of $(\overline{\gamma},\overline{\nu}_1,\overline{\nu}_2)$ is given by
\begin{align*}
  \overline{\ell}(t)&=\ell(t)\cos\theta(t)-n(t)\sin\theta(t), \ \overline{m}(t)=\dot\theta(t)+m(t), \ \overline{n}(t)=\ell(t)\sin\theta(t)+n(t)\cos\theta(t),\\
  \overline{\alpha}(t)&=(\alpha(t)+\dot\lambda(t)+\eta(t)m(t))\cos\theta(t)+(\lambda(t)m(t)-\dot\eta(t))\sin \theta(t). 
\end{align*}
\end{proposition}

\begin{theorem}\label{n1n2mu-Bertrand-type}
$(\gamma,\nu_1,\nu_2):I \to \R^3 \times \Delta$ is a $(\nu_1,\nu_2,\overline{\mu})$-Bertrand framed curve if and only if there exists a smooth map $(\lambda,\eta):I \to \R^2$ with $(\lambda, \eta) \not\equiv (0, 0)$ such that
\begin{align}\label{n1n2mu-Bertrand-type-condition}
\dot\lambda(t)-\eta(t)\ell(t)=0, \ \lambda(t)\ell(t)+\dot\eta(t)=0
\end{align}
for all $t \in I$. 
\end{theorem}
\demo
Suppose that $(\gamma,\nu_1,\nu_2):I \to \R^3 \times \Delta$ is a $(\nu_1,\nu_2,\overline{\mu})$-Bertrand framed curve. 
By differentiating $\overline{\gamma}(t)=\gamma(t)+\lambda(t)\nu_1(t)+\eta(t)\nu_2(t)$, we have 
$
\overline{\alpha}(t) \overline{\mu}(t)=(\dot\lambda(t)-\eta(t)\ell(t))\nu_1(t)+(\lambda(t)\ell(t)+\dot\eta(t))\nu_2(t)+(\alpha(t)+\lambda(t)m(t)+\eta(t)n(t))\mu(t)
$ for all $t \in I$. 
Since $\nu_1(t)\times\nu_2(t)=\mu(t)=\overline{\mu}(t)$, we have $\overline{\alpha}(t)=\alpha(t)+\lambda(t)m(t)+\eta(t)n(t)$, $\dot\lambda(t)-\eta(t)\ell(t)=0$ and $\lambda(t)\ell(t)+\dot\eta(t)=0$ for all $t \in I$. 
\par
Conversely, suppose that there exists a smooth map $(\lambda,\eta):I \to \R^2$ with $(\lambda, \eta) \not\equiv (0, 0)$ such that condition (\ref{n1n2mu-Bertrand-type-condition}) satisfies. 
Let $(\overline{\gamma},\overline{\nu}_1,\overline{\nu}_2):I \to \R^3 \times \Delta$ be 
$
\overline{\gamma}(t)=\gamma(t)+\lambda(t)\mu(t)+\eta(t)\nu_1(t), \ \overline{\nu}_1(t)=\cos\theta(t)\nu_1(t)-\sin\theta(t)\nu_2(t)
$ and 
$\overline{\nu}_2(t)=\sin\theta(t)\nu_1(t)+\cos\theta(t)\nu_2(t),
$
where $\theta:I\to\R$ is a smooth function. 
Since $\dot{\overline{\gamma}}(t)=(\alpha(t)+\lambda(t)m(t)+\eta(t)n(t))\mu(t)$, we have
$\dot{\overline{\gamma}}(t)\cdot\overline{\nu}_1(t)=\dot{\overline{\gamma}}(t)\cdot\overline{\nu}_2(t)=0$. 
It follows that $(\overline{\gamma},\overline{\nu}_1,\overline{\nu}_2)$ is a framed curve.
Moreover, by a direct calculation, we have $\overline{\mu}(t)=\overline{\nu}_1(t)\times\overline{\nu}_2(t)=(\cos^2\theta(t)+\sin^2\theta(t))\mu(t)=\mu(t)=\nu_1(t)\times\nu_2(t)$. 
Therefore, $(\gamma,\nu_1,\nu_2)$ is a $(\nu_1,\nu_2,\overline{\mu})$-Bertrand framed curve.
\enD

\begin{remark}\label{n1n2mu-Bertrand-type-equivalent}
$(\gamma,\nu_1,\nu_2)$ is a $(\nu_2,\nu_1,\overline{\mu})$-Bertrand framed curve if and only if $(\gamma,\nu_1,\nu_2)$ is a $(\nu_1,\nu_2,\overline{\mu})$-Bertrand framed curve.
\end{remark}
By the same method of Proposition \ref{n1n2n2-Bertrand-type_curvature}, we have the following.
\begin{proposition}\label{n1n2mu-Bertrand-type_curvature}
Suppose that $(\gamma,\nu_1,\nu_2)$ and $(\overline{\gamma},\overline{\nu}_1,\overline{\nu}_2):  I \to \R^3 \times \Delta$ are {$(\nu_1,\nu_2,\overline{\mu})$-mates}, where $\overline{\gamma}(t)=\gamma(t)+\lambda(t)\nu_1(t)+\eta(t)\nu_2(t), \overline{\nu}_1(t)=\cos\theta(t)\nu_1(t)-\sin\theta(t)\nu_2(t)$ and $\overline{\nu}_2(t)=\sin\theta(t)\nu_1(t)+\cos\theta(t)\nu_2(t)$. 
Then the curvature $(\overline{\ell}, \overline{m}, \overline{n}, \overline{\alpha})$ of $(\overline{\gamma},\overline{\nu}_1,\overline{\nu}_2)$ is given by
\begin{align*}
  \overline{\ell}(t)&=\ell(t)-\dot\theta(t), \ \overline{m}(t)=m(t)\cos\theta(t)-n(t)\sin\theta(t), \\
  \overline{n}(t)&=m(t)\sin\theta(t)+n(t)\cos\theta(t), \ \overline{\alpha}(t)=\alpha(t)+\lambda(t)m(t)+\eta(t)n(t). 
\end{align*}
\end{proposition}

\section{Evolutes and involutes of non-degenerate curves}
The evolutes are classical object and it is well-known properties of evolutes of regular curves (cf. \cite{Arnold1, Blaschke, Fuchs, Gray, Porteous, Romero-Sanabria, Uribe-Vargas1}).
Using $(\bv,\bw,\overline{\bx})$-Bertrand type curves, we define the involute of non-degenerate curves directly. 
Let $\gamma:I \to \R^3$ be a non-degenerate curve with the curvature $\kappa$ and the torsion $\tau$. 
\begin{definition}\label{involute_non-deg}{\rm
The involute ${I}nv(\gamma):I\rightarrow\R^3$ of non-degenerate curve $\gamma$ with $\tau\neq0$ is given by 
\begin{align*}
{I}nv(\gamma)(t)=\gamma(t)+\lambda^{I}(t)\bt(t)+\eta^{I}(t)\bn(t),
\end{align*}
where there exists a smooth map $(\lambda^{I}, \eta^{I}):I\to\R^2$ with $\eta^{I}\neq0$ 
such that 
\begin{align}\label{involute_non-deg-condition}
\begin{cases}
|\dot\gamma(t)|+\dot\lambda^{I}(t)-\eta^{I}(t)|\dot\gamma(t)|\kappa(t)=0,\\
\lambda^{I}(t)|\dot\gamma(t)|\kappa(t)+\dot\eta^{I}(t)=0
\end{cases}
\end{align}
for all $t \in I$.
}
\end{definition}
\begin{remark}{\rm
The involute ${I}nv(\gamma)$ and $\gamma$ are corresponding to the $(\bt,\bn,\overline{\bt})$-mates. 
}
\end{remark}

We consider conditions that evolutes and involutes are inverse operations of non-degenerate curves.
\begin{theorem}\label{evo_invo_non-deg_inverse}
Let $\gamma:I \to \R^3$ be a non-degenerate curve with curvature $\kappa$ and torsion $\tau$. 
Suppose that $\tau\neq0$. 
\par
$(1)$ Suppose that there exists a smooth map $(\lambda^{I}, \eta^{I}):I\to\R^2$ with $\eta^{I}\neq0$ such that condition \eqref{involute_non-deg-condition} holds. 
Then the evolute of the involute $Ev(Inv(\gamma))$ exists. Moreover, $Ev(Inv(\gamma))=\gamma$.
\par
$(2)$ Suppose that there exists a smooth map $(\lambda^{E}, \eta^{E}):I\to\R^2$ such that $\lambda^{E}=1/\kappa$ and $\eta^{E}=-{\dot{\kappa}}/({\vert \dot{\gamma} \vert \kappa^2 \tau})$, and equation \eqref{evolute_non-deg} holds. 
\par
$(i)$ Suppose that $\tau>0$ and $h>0$. 
If we take $(\lambda^{I}, \eta^{I}):I\to\R^2$ by $\lambda^{I}=-\eta^{E}$ and $\eta^{I}=\lambda^{E}$, then the involute of the evolute $Inv(Ev(\gamma))$ exists. 
Moreover, $Inv(Ev(\gamma))=\gamma$. 
\par
$(ii)$ Suppose that $\tau>0$ and $h<0$.
If we take $(\lambda^{I}, \eta^{I}):I\to\R^2$ by $\lambda^{I}=\eta^{E}$ and $\eta^{I}=-\lambda^{E}$, then the involute of the evolute $Inv(Ev(\gamma))$ exists. 
Moreover, $Inv(Ev(\gamma))=\gamma$. 
\par
$(iii)$ Suppose that $\tau<0$ and $h>0$. 
If we take $(\lambda^{I}, \eta^{I}):I\to\R^2$ by $\lambda^{I}=-\eta^{E}$ and $\eta^{I}=-\lambda^{E}$, then the involute of the evolute $Inv(Ev(\gamma))$ exists. 
Moreover, $Inv(Ev(\gamma))=\gamma$. 
\par
$(iv)$ Suppose that $\tau<0$ and $h<0$. 
If we take $(\lambda^{I}, \eta^{I}):I\to\R^2$ by $\lambda^{I}=\eta^{E}$ and $\eta^{I}=\lambda^{E}$, then the involute of the evolute $Inv(Ev(\gamma))$ exists. 
Moreover, $Inv(Ev(\gamma))=\gamma$. 
\end{theorem}
\demo
$(1)$ By Definition \ref{involute_non-deg} $(1)$, the involute $Inv(\gamma)$ is given by $Inv(\gamma)(t)=\gamma(t)+\lambda^{I}(t)\bt(t)+\eta^{I}(t)\bn(t)$, where $(\lambda^{I},\eta^{I}):I\to\R^2$ with $\eta^{I}\neq0$ is a smooth map such that condition (\ref{involute_non-deg-condition}) holds.  
By Remark \ref{tnt-Bertrand-type_curvature}, the curvature and the torsion of $Inv(\gamma)$ are given by $\kappa^{I}(t)=1/|\eta^{I}(t)|$ and $\tau^{I}(t)=\kappa(t)/(\eta^{I}(t)\tau(t))$.
By Theorem \ref{nbb-Bertrand-type} $(1)$, the evolute of the involute $Ev(Inv(\gamma))$ is given by 
\begin{align*}
Ev(Inv(\gamma))(t)&=Inv(\gamma)(t)+\lambda^{E}(t){\bn}^{I}(t)+\eta^{E}(t){\bb}^{I}(t).
\end{align*}
Since $\tau\neq0$ and $\eta^{I}\neq0$, we divide four cases.
\par
Suppose that $\tau>0$ and $\eta^{I}>0$.
Then the moving frame and the curvature of $Inv(\gamma)$ are given by $\{\bt^{I},\bn^{I},\bb^{I}\}=\{\bb,-\bn,\bt\}$ and $\kappa^{I}(t)=1/\eta^{I}(t)$. We also have $|\dot{Inv}(\gamma)(t)|=\eta^{I}(t)|\dot\gamma(t)|\tau(t)$. 
The evolute of the involute $Ev(Inv(\gamma))$ is given by 
$$
Ev(Inv(\gamma))(t)=\gamma(t)+(\lambda^{I}(t)+\eta^{E}(t))\bt(t)+(\eta^{I}(t)-\lambda^{E}(t))\bn(t).
$$
The condition (\ref{evolute_non-deg}) of $Inv(\gamma)$ is given by
\begin{align*}
h^{I}(t)&=\frac{|\dot{Inv}(\gamma)(t)|\tau^{I}(t)}{\kappa^{I}(t)}-\frac{d}{dt}\biggl(\frac{\dot\kappa^{I}(t)}{{\kappa^{I}}^2(t)}\frac{1}{|\dot{Inv}(\gamma)(t)|\tau^{I}(t)}\biggr)\\
&=\eta^{I}(t)|\dot\gamma(t)|\tau(t)\frac{\kappa(t)}{\eta^{I}(t)\tau(t)}\eta^{I}(t)-\frac{d}{dt}\biggl(-\frac{\dot\eta^{I}(t)}{(\eta^{I}(t))^2}(\eta^{I}(t))^2\frac{1}{\eta^{I}(t)|\dot\gamma(t)|\tau(t)}\frac{\eta^{I}(t)\tau(t)}{\kappa(t)}\biggr)\\
&=\eta^{I}(t)|\dot\gamma(t)|\kappa(t)-\frac{d}{dt}\biggl(-\frac{\dot\eta^{I}(t)}{|\dot\gamma(t)|\kappa(t)}\biggr)\\
&=\eta^{I}(t)|\dot\gamma(t)|\kappa(t)-\dot\lambda^{I}(t)\\
&=|\dot\gamma(t)|\neq0.
\end{align*} 
Thus, the evolute of the involute $Ev(Inv(\gamma))$ exists.
Moreover, by a direct calculation, we have 
\begin{align*}
\lambda^{E}(t)&=\frac{1}{\kappa^{I}(t)}=\eta^{I}(t), \ 
\eta^{E}(t)=\frac{\dot\lambda^{E}(t)}{|\dot{Inv}(\gamma)(t)|\tau^{I}(t)}=-\frac{\lambda^{I}(t)|\dot\gamma(t)|\kappa(t)}{|\dot\gamma(t)|\kappa(t)}=-\lambda^{I}(t).
\end{align*}
Therefore, we have $Ev(Inv(\gamma))(t)=\gamma(t)+(\lambda^{I}(t)+\eta^{E}(t))\bt(t)+(\eta^{I}(t)-\lambda^{E}(t))\bn(t)=\gamma(t)$.
\par
If $\tau>0$ and $\eta^{I}<0$ (respectively, $\tau<0$ and $\eta^{I}>0$, or $\tau<0$ and $\eta^{I}<0$), the moving frame and the curvature of $Inv(\gamma)$ are given by $\{\bt^{I},\bn^{I},\bb^{I}\}=\{-\bb,\bn,\bt\}$ (respectively, $\{-\bb,-\bn,-\bt\}$ or $\{\bb,\bn,-\bt\}$) and $\kappa^{I}(t)=-1/\eta^{I}(t)$ (respectively, $1/\eta^{I}(t)$ or $-1/\eta^{I}(t)$). We also have $|\dot{Inv}(\gamma)(t)|=-\eta^{I}(t)|\dot\gamma(t)|\tau(t)$ (respectively, $-\eta^{I}(t)|\dot\gamma(t)|\tau(t)$ or $\eta^{I}(t)|\dot\gamma(t)|\tau(t)$).
By the same calculations as in $\tau>0$ and $\eta^{I}>0$, we have $Ev(Inv(\gamma))(t)=\gamma(t)$.
\par
$(2)$ By Theorem \ref{nbb-Bertrand-type} $(1)$, the evolute $Ev(\gamma)$ is given by $Ev(\gamma)(t)=\gamma(t)+\lambda^{E}(t)\bn(t)+\eta^{E}(t)\bb(t)$, where $(\lambda^{E},\eta^{E}):I\to\R^2$ is a smooth map given by $\lambda^{E}(t)=1/\kappa(t)$ and $\eta^{E}(t)=-{\dot{\kappa}(t)}/({\vert \dot{\gamma}(t) \vert \kappa^2 (t)\tau(t)})$, and equation \eqref{evolute_non-deg} holds. 
By Remark \ref{nbb-Bertrand-type_curvature}, the curvature and the torsion of $Ev(\gamma)$ are given by $\kappa^{E}(t)={|\dot\gamma(t)||\tau(t)|}/{|h(t)|}$ and $\tau^{E}(t)={|\dot\gamma(t)|\kappa(t)}/{h(t)}$.
By Definition \ref{involute_non-deg}, the involute of the evolute $Inv(Ev(\gamma))$ is given by 
\begin{align*}
Inv(Ev(\gamma))(t)&=Ev(\gamma)(t)+\lambda^{I}(t){\bt}^{E}(t)+\eta^{I}(t){\bn}^{I}(t).
\end{align*}
Moreover, the condition (\ref{involute_non-deg-condition}) of $Ev(\gamma)$ is given by
\begin{align}\label{Inv(Ev)_non-deg}
\begin{cases}
|\dot{Ev}(\gamma)(t)|+\dot\lambda^{I}(t)-\eta^{I}(t)|\dot{Ev}(\gamma)(t)|\kappa^{E}(t)=0,\\
\lambda^{I}(t)|\dot{Ev}(\gamma)(t)|\kappa^{E}(t)+\dot\eta^{I}(t)=0.
\end{cases}
\end{align}
\par
$(i)$ Suppose that $\tau>0$ and $h>0$. 
Then the moving frame and the curvature of $Ev(\gamma)$ are given by $\{\bt^{E},\bn^{E},\bb^{E}\}=\{\bb,-\bn,\bt\}$ and $\kappa^{E}(t)=|\dot{\gamma}(t)|\tau(t)/h(t)$. We also have $|\dot{Ev}(\gamma)(t)|=h(t)$.  
The involute of the evolute $Inv(Ev(\gamma))$ is given by
$$
Inv(Ev(\gamma))(t)=\gamma(t)+(\lambda^{E}(t)-\eta^{I}(t))\bn(t)+(\eta^{E}(t)+\lambda^{I}(t))\bb(t).
$$
and a direct calculation, condition \eqref{Inv(Ev)_non-deg} is given by 
\begin{align*}
\begin{cases}
|\dot\gamma(t)|\tau(t)(\lambda^{E}(t)-\eta^{I}(t))+\dot{\eta}^{E}(t)+\dot\lambda^{I}(t)=0,\\
\lambda^{I}(t)|\dot\gamma(t)|\tau(t)+\dot\eta^{I}(t)=0.
\end{cases}
\end{align*}
If we take $\lambda^{I}(t)=-\eta^{E}(t)$ and $\eta^{I}(t)=\lambda^{E}(t)$, then condition (\ref{Inv(Ev)_non-deg}) satisfies. Thus, the involute of the evolute $Inv(Ev(\gamma))$ exists. 
Therefore, we have $Inv(Ev(\gamma))(t)=\gamma(t)+(\lambda^{E}(t)-\eta^{I}(t))\bn(t)+(\eta^{E}(t)+\lambda^{I}(t))\bb(t)=\gamma(t)$.
\par
$(ii)$, $(iii)$, $(iv)$ Suppose that $\tau>0$ and $h<0$ (respectively, $\tau<0$ and $h>0$, or $\tau<0$ and $h<0$). 
Then the moving frame and the curvature of $Ev(\gamma)$ are given by $\{\bt^{E},\bn^{E},\bb^{E}\}=\{-\bb,\bn,\bt\}$ (respectively, $\{\bb,\bn,-\bt\}$ or $\{-\bb,-\bn,-\bt\}$) and $\kappa^{E}(t)=-|\dot{\gamma}(t)|\tau(t)/h(t)$ (respectively, $-|\dot{\gamma}(t)|\tau(t)/h(t)$ or $|\dot{\gamma}(t)|\tau(t)/h(t)$).  
We also have $|\dot{Ev}(\gamma)(t)|=-h(t)$ (respectively, $h(t)$ or $-h(t)$). 
By the same calculations as in $(i)$, we have $Inv(Ev(\gamma))(t)=\gamma(t)$.
\enD
We consider the relation between involutes of non-degenerate curves.
\begin{lemma}\label{relation_Invs}
Let $\gamma:I \to \R^3$ be a non-degenerate curve with  the curvature $\kappa$ and the torsion $\tau$. 
Suppose that $\tau\neq0$ and there exist involutes $Inv_1(\gamma)$ and $Inv_2(\gamma)$, where $(\lambda^{I}_1, \eta^{I}_1):I\to\R^2$ and $(\lambda^{I}_2, \eta^{I}_2):I\to\R^2$ with $\eta^{I}_1\neq0$ and $\eta^{I}_2\neq0$ are smooth maps such that condition \eqref{involute_non-deg-condition} holds. 
Moreover, suppose that $(\lambda^{I}_1-\lambda^{I}_2, \eta^{I}_1-\eta^{I}_2) \not\equiv (0, 0)$.
Then the involutes $Inv_1(\gamma)$ and $Inv_2(\gamma)$ are $(\bn,\bb,\overline{\bt})$-mates.
\end{lemma}
\demo
Suppose that $\tau>0$. Then the moving frame of $Inv_1(\gamma)$ and $Inv_2(\gamma)$ are given by $\{\bt^{I}_1,\bn^{I}_1,\bb^{I}_1\}=\{\bt^{I}_2,\bn^{I}_2,\bb^{I}_2\}=\{\pm\bb,\mp\bn,\bt\}$.
By Definition \ref{Bertrand-type-regular2}, if we take $\lambda(t)=\pm(\eta^{I}_1(t)-\eta^{I}_2(t))$ and $\eta(t)=-(\lambda^{I}_1(t)-\lambda^{I}_2(t))$, then $(\lambda, \eta)\not\equiv(0,0)$. 
We denote $\overline{\gamma}(t)=Inv_1(\gamma)(t)+\lambda(t)\bn^{I}_1(t)+\eta(t)\bb^{I}_1(t)$. 
By a direct calculation, we have 
\begin{align*}
\overline{\gamma}(t)&=Inv_1(\gamma)(t)+\lambda(t)\bn_1^I(t)+\eta(t)\bb_1^I(t)\\
&=\gamma(t)+\lambda_1^{I}(t)\bt(t)+\eta_1^{I}(t)\bn(t)-(\eta^{I}_1(t)-\eta^{I}_2(t))\bn(t)-(\lambda^{I}_1(t)-\lambda^{I}_2(t))\bt(t)\\
&=\gamma(t)+\lambda_2^{I}(t)\bt(t)+\eta_2^{I}(t)\bn(t)\\
&=Inv_2(\gamma)(t).
\end{align*}
Therefore, $\overline{\gamma}$ is non-degenerate. 
Moreover, since $\overline{\bt}(t)=\bt^{I}_2(t)=\bt^{I}_1(t)=\bn^{I}_1(t)\times\bb^{I}_1(t)$, $\overline{\gamma}=Inv_2(\gamma)$ and $Inv_1(\gamma)$ are $(\bn,\bb,\overline{\bt})$-mates.
\par
On the other hand,  suppose that $\tau<0$. Then the moving frame $\{\bt^{I}_1,\bn^{I}_1,\bb^{I}_1\}=\{\bt^{I}_2,\bn^{I}_2,\bb^{I}_2\}=\{\pm\bb,\pm\bn,\bt\}$ and we obtain the same result as $\tau>0$.
\enD
We consider cases in which evolutes and involutes are not inverse operations. 
By using Theorem \ref{evo_invo_non-deg_inverse} $(2)$ and Lemma \ref{relation_Invs}, we have the following.
\begin{theorem}\label{Inv(Ev)_another}
Let $\gamma:I \to \R^3$ be a non-degenerate curve with  the curvature $\kappa$ and the torsion $\tau$. 
Suppose that $\tau\neq0$, equation \eqref{evolute_non-deg} holds and there exists a smooth map $(\lambda^{E}, \eta^{E}):I\to\R^2$ such that $\lambda^{E}=1/\kappa$ and $\eta^{E}=-{\dot{\kappa}}/({\vert \dot{\gamma} \vert \kappa^2 \tau})$. 
\par
$(i)$ Suppose that $\tau>0$ and $h>0$ and there exists a smooth map $(\lambda^{I},\eta^{I}):I\to\R^2$ with $(\lambda^{I}+\eta^{E},\eta^{I}-\lambda^{E})\not\equiv(0,0)$ such that condition \eqref{Inv(Ev)_non-deg} holds. 
Then $\gamma$ and $Inv(Ev(\gamma))$ are $(\bn,\bb,\overline{\bt})$-mates. 
\par
$(ii)$ Suppose that $\tau>0$ and $h<0$ and there exists a smooth map $(\lambda^{I},\eta^{I}):I\to\R^2$ with $(\lambda^{I}-\eta^{E},\eta^{I}+\lambda^{E})\not\equiv(0,0)$ such that condition \eqref{Inv(Ev)_non-deg} holds. 
Then $\gamma$ and $Inv(Ev(\gamma))$ are $(\bn,\bb,\overline{\bt})$-mates. 
\par
$(iii)$ Suppose that $\tau<0$ and $h>0$ and there exists a smooth map $(\lambda^{I},\eta^{I}):I\to\R^2$ with $(\lambda^{I}+\eta^{E},\eta^{I}+\lambda^{E})\not\equiv(0,0)$ such that condition \eqref{Inv(Ev)_non-deg} holds. 
Then $\gamma$ and $Inv(Ev(\gamma))$ are $(\bn,\bb,\overline{\bt})$-mates. 
\par
$(iv)$ Suppose that $\tau<0$ and $h<0$ and there exists a smooth map $(\lambda^{I},\eta^{I}):I\to\R^2$ with $(\lambda^{I}-\eta^{E},\eta^{I}-\lambda^{E})\not\equiv(0,0)$ such that condition \eqref{Inv(Ev)_non-deg} holds. 
Then $\gamma$ and $Inv(Ev(\gamma))$ are $(\bn,\bb,\overline{\bt})$-mates. 
\end{theorem}
\demo
$(i)$ By assumption, there exists the involute of the evolute $Inv(Ev(\gamma))(t)=\gamma(t)+(\lambda^{E}(t)-\eta^{I}(t))\bn(t)+(\eta^{E}(t)+\lambda^{I}(t))\bb(t)$, where $(\lambda^{I},\eta^{I}):I\to\R^2$ with $(\lambda^{I}+\eta^{E},\eta^{I}-\lambda^{E})\not\equiv(0,0)$ is a smooth map such that condition \eqref{Inv(Ev)_non-deg} holds. 
By Theorem \ref{evo_invo_non-deg_inverse} $(2)$, if we take $\lambda^{I}_2(t)=-\eta^{E}(t)$ and $\eta^{I}(t)=\lambda^{E}(t)$, then condition \eqref{Inv(Ev)_non-deg} also holds and we have $Inv_2(Ev(\gamma))(t)=\gamma(t)$. 
By Lemma \ref{relation_Invs}, $Inv(Ev(\gamma))$ and $Inv_2(Ev(\gamma))=\gamma$ are $(\bn,\bb,\overline{\bt})$-mates. 
\par
$(ii)$, $(iii)$, $(iv)$ We can also prove by the same method of $(i)$.
\enD

\section{Evolutes and involutes of framed curves}
Using $(\bv,\bw,\overline{\bx})$-Bertrand framed curves, we define evolutes and involutes of framed curves directly. 
\begin{definition}\label{evolute_involute_framed}{\rm
Let $(\gamma,\nu_1,\nu_2):I \to \R^3 \times \Delta$ be a framed curve with curvature $(\ell,m,n,\alpha)$. 
\par
$(1)$ The framed curve $\mathcal{E}(\gamma,\nu_1,\nu_2)=(\Ev(\gamma,\nu_1,\nu_2),\nu^{\mathcal{E}}_1(\gamma,\nu_1,\nu_2),\nu^{\mathcal{E}}_2(\gamma,\nu_1,\nu_2)):I \to \R^3 \times \Delta$ is given by 
\begin{align*}
&\Ev(\gamma,\nu_1,\nu_2)(t)=\gamma(t)+\lambda^{\mathcal{E}}(t)\nu_1(t)+\eta^{\mathcal{E}}(t)\nu_2(t),\\
&\nu^{\mathcal{E}}_1(\gamma,\nu_1,\nu_2)(t)=\sin\theta^{\mathcal{E}}(t)\nu_1(t)+\cos\theta^{\mathcal{E}}(t)\nu_2(t),\\
&\nu^{\mathcal{E}}_2(\gamma,\nu_1,\nu_2)(t)=\mu(t),
\end{align*}
where there exist a smooth map $(\lambda^{\mathcal{E}}, \eta^{\mathcal{E}}):I \to \R^2$ and a smooth function $\theta^{\mathcal{E}}:I\to\R$ such that 
\begin{align}\label{evolute_framed-condition}
\begin{cases}
\alpha(t)+\lambda^{\mathcal{E}}(t)m(t)+\eta^{\mathcal{E}}(t)n(t)=0,\\
(\dot\lambda^{\mathcal{E}}(t)-\eta^{\mathcal{E}}(t)\ell(t))\sin\theta^{\mathcal{E}}(t)+(\lambda^{\mathcal{E}}(t)\ell(t)+\dot\eta^{\mathcal{E}}(t))\cos\theta^{\mathcal{E}}(t)=0
\end{cases}
\end{align}
for all $t \in I$. 
Then we say that $\Ev(\gamma,\nu_1,\nu_2):I \to \R^3$ is an {\it evolute} of the framed curve $(\gamma,\nu_1,\nu_2)$. 
\par
$(2)$ The framed curve $\mathcal{I}(\gamma,\nu_1,\nu_2)=(\Inv(\gamma,\nu_1,\nu_2),\nu^{\mathcal{I}}_1(\gamma,\nu_1,\nu_2),\nu^{\mathcal{I}}_2(\gamma,\nu_1,\nu_2)):I \to \R^3 \times \Delta$ is given by 
\begin{align*}
&\Inv(\gamma,\nu_1,\nu_2)(t)=\gamma(t)+\lambda^{\mathcal{I}}(t)\mu(t)+\eta^{\mathcal{I}}(t)\nu_1(t),\\
&\nu^{\mathcal{I}}_1(\gamma,\nu_1,\nu_2)(t)=\cos\theta^{\mathcal{I}}(t)\mu(t)-\sin\theta^{\mathcal{I}}(t)\nu_1(t),\\
&\nu^{\mathcal{I}}_2(\gamma,\nu_1,\nu_2)(t)=\sin\theta^{\mathcal{I}}(t)\mu(t)+\cos\theta^{\mathcal{I}}(t)\nu_1(t),
\end{align*}
where there exists a smooth map $(\lambda^{\mathcal{I}}, \eta^{\mathcal{I}}):I \to \R^2$ such that 
\begin{align}\label{involute_framed-condition}
\begin{cases}
-\lambda^{\mathcal{I}}(t)m(t)+\dot\eta^{\mathcal{I}}(t)=0,\\
\alpha(t)+\dot\lambda^{\mathcal{I}}(t)+\eta^{\mathcal{I}}(t)m(t)=0
\end{cases}
\end{align}
for all $t \in I$ and $\theta^{\mathcal{I}}:I\to\R$ is a smooth function. 
Then we say that $\Inv(\gamma,\nu_1,\nu_2):I \to \R^3$ is an {\it involute} of the framed curve $(\gamma,\nu_1,\nu_2)$. 
}
\end{definition}
\begin{remark}{\rm 
The framed curves $\mathcal{E}(\gamma,\nu_1,\nu_2)$ (respectively, $\mathcal{I}(\gamma,\nu_1,\nu_2)$) and $(\gamma,\nu_1,\nu_2)$ are corresponding to the $(\nu_1,\nu_2,\overline{\nu}_2)$ (respectively, $(\mu,\nu_1,\overline{\mu})$)-mates. 
}
\end{remark}

We consider the relation between the evolute $\Ev(\gamma,\nu_1,\nu_2)$ of a framed curve and circular evolute with respect to Bishop directions of a framed curve (cf. \cite{Honda-Takahashi-Preprint}).
For a framed curve $(\gamma,\nu_1,\nu_2):I \to \R^3 \times \Delta$, we call a unit vector $\bv(t)\in {\langle \nu_1(t),\nu_2(t) \rangle}_{\R}$ is a {\it Bishop direction} (or, {\it Bishop vector}) if there exists a smooth function $\beta:I\to\R$ such that $\dot\bv(t)=\beta(t)\mu(t)$. 
We also call a moving frame $\{\bv,\bw,\mu\}$ a {\it Bishop frame} along $\gamma$ is both $\bv$ and $\bw$ are Bishop directions. 
We define $(\bv,\bw):I\to\Delta$ by
$$
\begin{pmatrix}
{\bv}(t) \\
{\bw}(t)
\end{pmatrix}
=
\begin{pmatrix}
\cos \theta(t) & -\sin \theta(t)\\
\sin \theta(t) & \cos \theta(t)
\end{pmatrix}
\begin{pmatrix}
{\nu_1}(t) \\
{\nu_2}(t)
\end{pmatrix},
$$
where $\theta:I\to\R$ is a smooth function. Then $(\gamma,\bv,\bw):I \to \R^3 \times \Delta$ is also a framed curve with the unit tangent vector $\mu$. 
The moving frame $\{\bv,\bw,\mu\}$ is called a {\it rotated frame} of $\{\nu_1,\nu_2,\mu\}$ by $\theta$.
Then the Frenet-Serret type formula of the rotated frame $\{\bv,\bw,\mu\}$ is given by 
\begin{align}\label{Bishop frame}
\left(
\begin{array}{c}
\dot{\bv}(t)\\
\dot{\bw}(t)\\
\dot{\mu}(t)
\end{array} \right)=
\left(
\begin{array}{ccc}
0 & \overline{\ell}(t) & \overline{m}(t)\\
-\overline{\ell}(t) & 0 & \overline{n}(t)\\
-\overline{m}(t) & -\overline{n}(t) & 0
\end{array}\right)
\left(
\begin{array}{c}
\bv(t)\\
\bw(t)\\
\mu(t)
\end{array}\right), 
\ \dot{\gamma}(t)={\alpha}(t)\mu(t),
\end{align}
where $\overline{\ell}(t)=\ell(t)-\dot\theta(t)$, $\overline{m}(t)=m(t)\cos\theta(t)-n(t)\sin\theta(t)$ and $\overline{n}(t)=m(t)\sin\theta(t)+n(t)\cos\theta(t)$.
If we take a function $\theta:I\to\R$ which satisfies $\dot\theta(t)=\ell(t)$ for all $t\in I$, then the rotated frame $\{\bv,\bw,\mu\}$ is a Bishop frame.
\par
Let $(\gamma,\nu_1,\nu_2):I \to \R^3 \times \Delta$ be a framed curve with a Bishop frame $\{\bv,\bw,\mu\}$, namely, $\overline{\ell}(t)=0$ for all $t\in I$ in the Frenet-Serret type formula in \eqref{Bishop frame}.
\begin{definition}{\rm
We assume that $\overline{m}(t)\neq0$ for all $t\in I$. Then the ({\it circular}) {\it evolute} $E_\gamma[\bv]:I\to\R^3$ of $\gamma$ with respect to $\bv$ (or, {\it $\bv$-evolute}) is given by 
$$
E_\gamma[\bv](t)=\gamma(t)-\frac{\alpha(t)}{\overline{m}(t)}\bv(t).
$$
}
\end{definition}
Then $(E_\gamma[\bv],\bw,\mu):I\to\R^3\times\Delta$ is a framed curve with the curvature $(\overline{n},0,-\overline{m},-(d/dt)(\alpha/\overline{m}))$ and the moving frame along the circular evolute $E_\gamma[\bv]$ is given by $\{\bw,\mu,\bv\}$.
\begin{proposition}\label{evolute-evolute}
Let $(\gamma,\nu_1,\nu_2):I \to \R^3 \times \Delta$ be a framed curve with a Bishop frame $\{\bv,\bw,\mu\}$ and $\overline{m}(t)\neq0$ for all $t\in I$.
If we take a smooth map $(\lambda^{\mathcal{E}}, \eta^{\mathcal{E}}):I \to \R^2$ and a smooth function $\theta^{\mathcal{E}}:I\to\R$ such that $\lambda^{\mathcal{E}}=-(\alpha\cos\theta)/\overline{m}$, $\eta^{\mathcal{E}}=(\alpha\sin\theta)/\overline{m}$ and $\theta^{\mathcal{E}}=\theta$, then we have the framed curve $\mathcal{E}(\gamma,\nu_1,\nu_2)=(E_\gamma[\bv],\bw,\mu)$, namely, in the case, $\Ev(\gamma,\nu_1,\nu_2)$ is a circular evolute.
\end{proposition}
\demo
If we take a smooth map $(\lambda^{\mathcal{E}}, \eta^{\mathcal{E}}):I \to \R^2$ and a smooth function $\theta^{\mathcal{E}}:I\to\R$ such that
\begin{align*}
\lambda^{\mathcal{E}}(t)&=-\frac{\alpha(t)\cos\theta(t)}{\overline{m}(t)}=-\frac{\alpha(t)\cos\theta(t)}{m(t)\cos\theta(t)-n(t)\sin\theta(t)}, \\
\eta^{\mathcal{E}}(t)&=\frac{\alpha(t)\sin\theta(t)}{\overline{m}(t)}=\frac{\alpha(t)\sin\theta(t)}{m(t)\cos\theta(t)-n(t)\sin\theta(t)}
\end{align*}
and $\theta^{\mathcal{E}}(t)=\theta(t)$, then condition \eqref{evolute_framed-condition} satisfies. 
Thus, the evolute $\Ev(\gamma,\nu_1,\nu_2)$ of the framed curve $(\gamma,\nu_1,\nu_2)$ exists.
By Definition \ref{evolute_involute_framed} $(1)$, the framed curve 
$$
\mathcal{E}(\gamma,\nu_1,\nu_2)=(\Ev(\gamma,\nu_1,\nu_2),\nu^{\mathcal{E}}_1(\gamma,\nu_1,\nu_2),\nu^{\mathcal{E}}_2(\gamma,\nu_1,\nu_2))
$$
is given by 
\begin{align*}
\Ev(\gamma,\nu_1,\nu_2)(t)&=\gamma(t)+\lambda^{\mathcal{E}}(t)\nu_1(t)+\eta^{\mathcal{E}}(t)\nu_2(t)\\
&=\gamma(t)-\frac{\alpha(t)\cos\theta(t)}{\overline{m}(t)}\nu_1(t)+\frac{\alpha(t)\sin\theta(t)}{\overline{m}(t)}\nu_2(t)\\
&=\gamma(t)-\frac{\alpha(t)}{\overline{m}(t)}\left(\cos\theta(t)\nu_1(t)-\sin\theta(t)\nu_2(t)\right)\\
&=\gamma(t)-\frac{\alpha(t)}{\overline{m}(t)}\bv(t)=E_\gamma[\bv](t),\\
\nu^{\mathcal{E}}_1(\gamma,\nu_1,\nu_2)(t)&=\sin\theta^{\mathcal{E}}(t)\nu_1(t)+\cos\theta^{\mathcal{E}}(t)\nu_2(t)\\
&=\sin\theta(t)\nu_1(t)+\cos\theta(t)\nu_2(t)=\bw(t),\\
\nu^{\mathcal{E}}_2(\gamma,\nu_1,\nu_2)(t)&=\mu(t).
\end{align*}
It follows that we have $\mathcal{E}(\gamma,\nu_1,\nu_2)=(E_\gamma[\bv],\bw,\mu)$.
\enD
We consider the relation between the involute $\Inv(\gamma,\nu_1,\nu_2)$ of a framed curve and the involute of $\gamma$ respect to $t_0$ (or $t_0$-involute) (cf. \cite{Honda-Takahashi-Preprint}).
\begin{definition}{\rm
Let $(\gamma,\nu_1,\nu_2):I \to \R^3 \times \Delta$ be a framed curve with $m^2(t)+n^2(t)\neq0$ for all $t\in I$. 
Then the {\it involute} $I_\gamma[t_0]:I\to\R^3$ of $\gamma$ with respect to $t_0$ (or, {\it $t_0$-involute}) is given by 
$$
I_\gamma[t_0](t)=\gamma(t)-\left(\int_{t_0}^{t}\alpha(t)dt\right)\mu(t).
$$
}
\end{definition}
We define smooth maps $\bxi$ and $\beeta:I\to S^{2}$ by
$$
\bxi(t)=\frac{n(t)\nu_1(t)-m(t)\nu_2(t)}{\sqrt{m^2(t)+n^2(t)}}, \ \beeta(t)=\bxi(t)\times\mu(t)=\frac{-m(t)\nu_1(t)-n(t)\nu_2(t)}{\sqrt{m^2(t)+n^2(t)}}.
$$
Then $(I_\gamma[t_0],\bxi,\mu):I \to \R^3 \times \Delta$ is a framed curve with the curvature $(0,f,\sqrt{m^2+n^2},\\-(\int_{t_0}^{t}\alpha(t)dt)\sqrt{m^2+n^2})$, where 
$$
f(t)=\frac{\dot{m}(t)n(t)-m(t)\dot{n}(t)-\ell(t)(m^2(t)+n^2(t))}{\sqrt{m^2(t)+n^2(t)}}
$$
and the moving frame along the involute of $\gamma$ respect to $t_0$ is given by $\{\bxi,\mu,\beeta\}$.
Note that $\{\bxi,\mu,\beeta\}$ is one of the Bishop frames along $t_0$ .
\begin{proposition}\label{involute-involute}
Let $(\gamma,\nu_1,\nu_2):I \to \R^3 \times \Delta$ be a framed curve with $m^2(t)+n^2(t)\neq0$ for all $t\in I$ and $t_0\in I$. 
Suppose that there exists the framed curve $\mathcal{I}(\gamma,\nu_1,\nu_2):I \to \R^3 \times \Delta$, where $(\lambda^{\mathcal{I}}, \eta^{\mathcal{I}}):I \to \R^2$ is a smooth map such that condition $(\ref{involute_framed-condition})$ holds and $\theta^{\mathcal{I}}:I\to\R$ is a smooth function. 
If we take a smooth map $(\lambda^{\mathcal{I}}, \eta^{\mathcal{I}}):I \to \R^2$ such that $\lambda^{\mathcal{I}}\not\equiv0$ and $\eta^{\mathcal{I}}=0$, then we have the framed curve $\mathcal{I}(\gamma,\nu_1,\nu_2)=(I_\gamma[t_0],\bxi,\mu)$, namely, in the case, $\Inv(\gamma,\nu_1,\nu_2)$ is an involute of $\gamma$ respect to $t_0$.
\end{proposition}
\demo
If we take a smooth map $(\lambda^{\mathcal{I}}, \eta^{\mathcal{I}}):I \to \R^2$ such that $\lambda^{\mathcal{I}}(t)\not\equiv0$ and $\eta^{\mathcal{I}}(t)=0$, then condition $(\ref{involute_framed-condition})$ becomes $-\lambda^{\mathcal{I}}(t)m(t)=0$ and $\alpha(t)+\dot\lambda(t)=0$.
Thus, we have $\lambda^{\mathcal{I}}(t)=-\int_{t_0}^{t}\alpha(t)dt$ for $t_0\in I$ and $m(t)=0$ for all $t\in I$.
Since $n(t)\neq0$ for all $t\in I$, firstly, we suppose that $n(t)<0$ for all $t\in I$. 
Then the Bishop frame of $(I_\gamma[t_0],\bxi,\mu)$ is given by $\{\bxi,\mu,\beeta\}=\{-\nu_1,\mu,\nu_2\}$.
If we take $\theta^{\mathcal{I}}:I\to\R$ such that $\theta^{\mathcal{I}}(t)=\pi/2$, then we have $\nu_1^{\mathcal{I}}(\gamma,\nu_1,\nu_2)(t)=-\nu_1(t)$ and $\nu_2^{\mathcal{I}}(\gamma,\nu_1,\nu_2)(t)=\mu(t)$. 
Therefore, then $\mathcal{I}(\gamma,\nu_1,\nu_2)(t)=(\gamma(t)-(\int_{t_0}^{t}\alpha(t)dt)\mu(t),-\nu_1(t),\mu(t))=(I_\gamma[t_0],\bxi,\mu)(t)$. 
\par
The next, we suppose that $n(t)>0$ for all $t\in I$. 
Then the Bishop frame of $(I_\gamma[t_0],\bxi,\mu)$ is given by $\{\bxi,\mu,\beeta\}=\{\nu_1,\mu,-\nu_2\}$.
By Proposition \ref{frame-change} $(3)$ and $(\gamma,\nu_1,\nu_2)$ is a $(\mu,\nu_1,\overline{\mu})$-Bertrand framed curve, $(\gamma,\nu_1,\nu_2)$ is also a $(\mu,-\nu_1,-\overline{\mu})$-Bertrand framed curve.
Thus, we have 
\begin{align*}
&\mathcal{I}(\gamma,\nu_1,\nu_2)(t)\\
&=(\gamma(t)+\lambda^{\mathcal{I}}(t)\mu(t)-\eta^{\mathcal{I}}(t)\nu_1(t),\cos\theta^{\mathcal{I}}(t)\mu(t)-\sin\theta^{\mathcal{I}}(t)\nu_1(t),-\sin\theta^{\mathcal{I}}(t)\mu(t)-\cos\theta^{\mathcal{I}}(t)\nu_1(t)).
\end{align*}
If we take $\theta^{\mathcal{I}}:I\to\R$ such that $\theta^{\mathcal{I}}(t)=-\pi/2$, then we have $\nu_1^{\mathcal{I}}(\gamma,\nu_1,\nu_2)(t)=\nu_1(t)$ and $\nu_2^{\mathcal{I}}(\gamma,\nu_1,\nu_2)(t)=\mu(t)$. 
Therefore, then $\mathcal{I}(\gamma,\nu_1,\nu_2)(t)=(\gamma(t)-(\int_{t_0}^{t}\alpha(t)dt)\mu(t),\nu_1(t),\mu(t))=(I_\gamma[t_0],\bxi,\mu)(t)$. 
\enD
We consider conditions that evolutes and involutes are inverse operations of framed curves.

\begin{theorem}\label{evo_invo_framed_inverse}
Let $(\gamma,\nu_1,\nu_2):I \to \R^3 \times \Delta$ be a framed curve with curvature $(\ell,m,n,\alpha)$.\par
$(1)$ Suppose that there exists the framed curve $\mathcal{I}(\gamma,\nu_1,\nu_2):I \to \R^3 \times \Delta$, where $(\lambda^{\mathcal{I}}, \eta^{\mathcal{I}}):I \to \R^2$ is a smooth map such that condition $(\ref{involute_framed-condition})$ holds and $\theta^{\mathcal{I}}:I\to\R$ is a smooth function. 
If we take
$
\lambda^{\mathcal{E}}=\lambda^{\mathcal{I}}\cos\theta^{\mathcal{I}}-\eta^{\mathcal{I}}\sin\theta^{\mathcal{I}},\ \eta^{\mathcal{E}}=\lambda^{\mathcal{I}}\sin\theta^{\mathcal{I}}+\eta^{\mathcal{I}}\cos\theta^{\mathcal{I}}
$
and $\theta^{\mathcal{E}}=-\theta^{\mathcal{I}}$, then the evolute of the involute $\Ev(\mathcal{I}(\gamma,\nu_1,\nu_2))$ exists. Moreover, we have $\mathcal{E}(\mathcal{I}(\gamma,\nu_1,\nu_2))=(\gamma,\nu_1,\nu_2)$.
\par
$(2)$ Suppose that there exists the framed curve $\mathcal{E}(\gamma,\nu_1,\nu_2):I \to \R^3 \times \Delta$, where $(\lambda^{\mathcal{E}}, \eta^{\mathcal{E}}):I \to \R^2$ is a smooth map and  $\theta^{\mathcal{E}}:I\to\R$ is a smooth function such that condition $(\ref{evolute_framed-condition})$ holds. 
If we take 
$
\lambda^{\mathcal{I}}=\lambda^{\mathcal{E}}\cos\theta^{\mathcal{E}}-\eta^{\mathcal{E}}\sin\theta^{\mathcal{E}},\ \eta^{\mathcal{I}}=\lambda^{\mathcal{E}}\sin\theta^{\mathcal{E}}+\eta^{\mathcal{E}}\cos\theta^{\mathcal{E}},
$
then the involute of the evolute $\Inv(\mathcal{E}(\gamma,\nu_1,\nu_2))$ exists. Moreover, if we take $\theta^{\mathcal{I}}=-\theta^{\mathcal{E}}$, then we have $\mathcal{I}(\mathcal{E}(\gamma,\nu_1,\nu_2))=(\gamma,\nu_1,\nu_2)$.
\end{theorem}
\demo
$(1)$ By Definition \ref{evolute_involute_framed} (2), the framed curve $\mathcal{I}(\gamma,\nu_1,\nu_2)$ is given by 
\begin{align*}
&\mathcal{I}(\gamma,\nu_1,\nu_2)(t)\\
&=(\gamma(t)+\lambda^{\mathcal{I}}(t)\mu(t)+\eta^{\mathcal{I}}(t)\nu_1(t),\cos\theta^{\mathcal{I}}(t)\mu(t)-\sin\theta^{\mathcal{I}}(t)\nu_1(t),\sin\theta^{\mathcal{I}}(t)\mu(t)+\cos\theta^{\mathcal{I}}(t)\nu_1(t)),
\end{align*}
where $(\lambda^{\mathcal{I}}, \eta^{\mathcal{I}}):I \to \R^2$ is a smooth map such that condition $(\ref{involute_framed-condition})$ holds and $\theta^{\mathcal{I}}:I\to\R$ is a smooth function.
The condition (\ref{evolute_framed-condition}) of $\mathcal{I}(\gamma,\nu_1,\nu_2)$ is given by
\begin{align}\label{Ev(Inv)_framed}
\begin{cases}
(\eta^{\mathcal{I}}(t)-\lambda^{\mathcal{E}}(t)\sin\theta^{\mathcal{I}}(t)+\eta^{\mathcal{E}}(t)\cos\theta^{\mathcal{I}}(t))\ell(t)\\
\hspace{5mm}-(\lambda^{\mathcal{I}}(t)+\lambda^{\mathcal{E}}(t)\cos\theta^{\mathcal{I}}(t)+\eta^{\mathcal{E}}(t)\sin\theta^{\mathcal{I}}(t))n(t)=0,\\
\left(\dot\lambda^{\mathcal{E}}(t)+\eta^{\mathcal{E}}(t)(\dot\theta^{\mathcal{I}}(t)+m(t))\right)\sin\theta^{\mathcal{E}}(t)\\
\hspace{5mm}+\left(-\lambda^{\mathcal{E}}(t)(\dot\theta^{\mathcal{I}}(t)+m(t))+\dot\eta^{\mathcal{E}}(t)\right)\cos\theta^{\mathcal{E}}(t)=0.
\end{cases}
\end{align}
If we take $\lambda^{\mathcal{E}}(t)=\lambda^{\mathcal{I}}(t)\cos\theta^{\mathcal{I}}(t)-\eta^{\mathcal{I}}(t)\sin\theta^{\mathcal{I}}(t)$, $\eta^{\mathcal{E}}(t)=\lambda^{\mathcal{I}}(t)\sin\theta^{\mathcal{I}}(t)+\eta^{\mathcal{I}}(t)\cos\theta^{\mathcal{I}}(t)$ and $\theta^{\mathcal{E}}(t)=-\theta^{\mathcal{I}}(t)$, then condition (\ref{Ev(Inv)_framed}) satisfies. Thus, the evolute of the involute $\Ev(\mathcal{I}(\gamma,\nu_1,\nu_2))$ exists. 
By Definition \ref{evolute_involute_framed} (1), the framed curve $\mathcal{E}(\gamma,\nu_1,\nu_2)$ of the framed curve $\mathcal{I}(\gamma,\nu_1,\nu_2)$,
$$
\mathcal{E}(\mathcal{I}(\gamma,\nu_1,\nu_2))=\left(\Ev(\mathcal{I}(\gamma,\nu_1,\nu_2)),\nu_1^{\mathcal{E}}(\mathcal{I}(\gamma,\nu_1,\nu_2)),\nu_2^{\mathcal{E}}(\mathcal{I}(\gamma,\nu_1,\nu_2))\right)
$$
is given by 
\begin{align*}
\Ev(\mathcal{I}(\gamma,\nu_1,\nu_2))(t)&=\Inv(\gamma,\nu_1,\nu_2)(t)+\lambda^{\mathcal{E}}(t)\nu_1^{\mathcal{I}}(\gamma,\nu_1,\nu_2)(t)+\eta^{\mathcal{E}}(t)\nu_2^{\mathcal{I}}(\gamma,\nu_1,\nu_2)(t)\\
&=\gamma(t)+(\lambda^{\mathcal{I}}(t)+\lambda^{\mathcal{E}}(t)\cos\theta^{\mathcal{I}}(t)+\eta^{\mathcal{E}}(t)\sin\theta^{\mathcal{I}}(t))\mu(t)\\
&\hspace{5mm}+(\eta^{\mathcal{I}}(t)-\lambda^{\mathcal{E}}(t)\sin\theta^{\mathcal{I}}(t)+\eta^{\mathcal{E}}(t)\cos\theta^{\mathcal{I}}(t))\nu_1(t)\\
&=\gamma(t),\\
\nu_1^{\mathcal{E}}(\mathcal{I}(\gamma,\nu_1,\nu_2))(t)&=\sin\theta^{\mathcal{E}}(t)\nu_1^{\mathcal{I}}(\gamma,\nu_1,\nu_2)(t)+\cos\theta^{\mathcal{E}}(t)\nu_2^{\mathcal{I}}(\gamma,\nu_1,\nu_2)(t)\\
&=\sin(\theta^{\mathcal{E}}(t)+\theta^{\mathcal{I}}(t))\mu(t)+\cos(\theta^{\mathcal{E}}(t)+\theta^{\mathcal{I}}(t))\nu_1(t)=\nu_1(t),\\
\nu_2^{\mathcal{E}}(\mathcal{I}(\gamma,\nu_1,\nu_2))(t)&=\mu^{\mathcal{I}}(t)=\nu_2(t).
\end{align*}
\par
$(2)$ By Definition \ref{evolute_involute_framed} (1), the framed curve $\mathcal{E}(\gamma,\nu_1,\nu_2)$ is given by 
\begin{align*}
&\mathcal{E}(\gamma,\nu_1,\nu_2)(t)\\
&=(\gamma(t)+\lambda^{\mathcal{E}}(t)\nu_1(t)+\eta^{\mathcal{E}}(t)\nu_2(t),\sin\theta^{\mathcal{E}}(t)\nu_1(t)+\cos\theta^{\mathcal{E}}(t)\nu_2(t),\mu(t)),
\end{align*}
where $(\lambda^{\mathcal{I}}, \eta^{\mathcal{I}}):I \to \R^2$ is a smooth map and $\theta^{\mathcal{E}}:I\to\R$ is a smooth function such that condition $(\ref{evolute_framed-condition})$ holds.
The condition (\ref{involute_framed-condition}) of $\mathcal{E}(\gamma,\nu_1,\nu_2)$ is given by
\begin{align}\label{Inv(Ev)_framed}
\begin{cases}
-\lambda^{\mathcal{I}}(t)(\dot\theta^{\mathcal{E}}(t)-\ell(t))+\dot\eta^{\mathcal{I}}(t)=0,\\
\dot\lambda^{\mathcal{I}}(t)+\eta^{\mathcal{I}}(t)(\dot\theta^{\mathcal{E}}(t)-\ell(t))\\
\hspace{5mm}+(\dot\lambda^{\mathcal{E}}(t)-\eta^{\mathcal{E}}(t)\ell(t))\cos\theta^{\mathcal{E}}(t)-(\lambda^{\mathcal{E}}(t)\ell(t)+\dot\eta^{\mathcal{E}}(t))\sin\theta^{\mathcal{E}}(t)=0.
\end{cases}
\end{align}
If we take $\lambda^{\mathcal{I}}(t)=\lambda^{\mathcal{E}}(t)\cos\theta^{\mathcal{E}}(t)-\eta^{\mathcal{E}}(t)\sin\theta^{\mathcal{E}}(t)$ and $\eta^{\mathcal{I}}(t)=\lambda^{\mathcal{E}}(t)\sin\theta^{\mathcal{E}}(t)+\eta^{\mathcal{E}}(t)\cos\theta^{\mathcal{E}}(t)$, then condition (\ref{Inv(Ev)_framed}) satisfies. Thus, the involute of the evolute $\Inv(\mathcal{E}(\gamma,\nu_1,\nu_2))$ exists. 
By Definition \ref{evolute_involute_framed} (2), the framed curve $\mathcal{I}(\gamma,\nu_1,\nu_2)$ of the framed curve $\mathcal{E}(\gamma,\nu_1,\nu_2)$,
$$\mathcal{I}(\mathcal{E}(\gamma,\nu_1,\nu_2))=\left(\Inv(\mathcal{E}(\gamma,\nu_1,\nu_2)),\nu_1^{\mathcal{I}}(\mathcal{E}(\gamma,\nu_1,\nu_2)),\nu_2^{\mathcal{I}}(\mathcal{E}(\gamma,\nu_1,\nu_2))\right)$$
is given by 
\begin{align*}
\Inv(\mathcal{E}(\gamma,\nu_1,\nu_2))(t)&=\Ev(\gamma,\nu_1,\nu_2)(t)+\lambda^{\mathcal{I}}(t)\mu^{\mathcal{E}}(\gamma,\nu_1,\nu_2)(t)+\eta^{\mathcal{I}}(t)\nu_1^{\mathcal{E}}(\gamma,\nu_1,\nu_2)(t)\\
&=\gamma(t)+(\lambda^{\mathcal{E}}(t)+\lambda^{\mathcal{I}}(t)\cos\theta^{\mathcal{E}}(t)+\eta^{\mathcal{I}}(t)\sin\theta^{\mathcal{E}}(t))\nu_1(t)\\
&\hspace{5mm}+(\eta^{\mathcal{E}}(t)-\lambda^{\mathcal{I}}(t)\sin\theta^{\mathcal{E}}(t)+\eta^{\mathcal{I}}(t)\cos\theta^{\mathcal{E}}(t))\nu_2(t)\\
&=\gamma(t),\\
\nu_1^{\mathcal{I}}(\mathcal{E}(\gamma,\nu_1,\nu_2))(t)&=\cos\theta^{\mathcal{I}}(t)\mu^{\mathcal{E}}(\gamma,\nu_1,\nu_2)(t)-\sin\theta^{\mathcal{I}}(t)\nu_1^{\mathcal{E}}(\gamma,\nu_1,\nu_2)(t)\\
&=\cos(\theta^{\mathcal{E}}(t)+\theta^{\mathcal{I}}(t))\nu_1(t)-\sin(\theta^{\mathcal{E}}(t)+\theta^{\mathcal{I}}(t))\nu_2(t),\\
\nu_2^{\mathcal{I}}(\mathcal{E}(\gamma,\nu_1,\nu_2))(t)&=\sin\theta^{\mathcal{I}}(t)\mu^{\mathcal{E}}(\gamma,\nu_1,\nu_2)(t)+\cos\theta^{\mathcal{I}}(t)\nu_1^{\mathcal{E}}(\gamma,\nu_1,\nu_2)(t)\\
&=\sin(\theta^{\mathcal{E}}(t)+\theta^{\mathcal{I}}(t))\nu_1(t)+\cos(\theta^{\mathcal{E}}(t)+\theta^{\mathcal{I}}(t))\nu_2(t).
\end{align*}
Moreover, if we take $\theta^{\mathcal{I}}(t)=-\theta^{\mathcal{E}}(t)$, then we have $\mathcal{I}(\mathcal{E}(\gamma,\nu_1,\nu_2))(t)=(\gamma,\nu_1,\nu_2)(t)$.
\enD
We consider the relation between the evolute and evolute (respectively, involute and involute) of framed curves.
\begin{lemma}\label{relation_Invs_Evs}
Let $(\gamma,\nu_1,\nu_2):I \to \R^3 \times \Delta$ be a framed curve with curvature $(\ell,m,n,\alpha)$.\par
$(1)$ Suppose that there exist the framed curves $\mathcal{E}_1(\gamma,\nu_1,\nu_2), \mathcal{E}_2(\gamma,\nu_1,\nu_2):I \to \R^3 \times \Delta$, where $\theta^{\mathcal{E}}_1, \theta^{\mathcal{E}}_2:I\to\R$ are smooth functions and $(\lambda^{\mathcal{E}}_1, \eta^{\mathcal{E}}_1), (\lambda^{\mathcal{E}}_2, \eta^{\mathcal{E}}_2):I\to\R^2$ are smooth maps such that condition \eqref{evolute_framed-condition} holds. 
Moreover, suppose that $(\lambda^{\mathcal{E}}_1-\lambda^{\mathcal{E}}_2, \eta^{\mathcal{E}}_1-\eta^{\mathcal{E}}_2) \not\equiv (0, 0)$.
Then the framed curves $\mathcal{E}_1(\gamma,\nu_1,\nu_2)$ and $\mathcal{E}_2(\gamma,\nu_1,\nu_2)$ are $(\mu,\nu_1,\overline{\nu}_2)$-mates.
\par
$(2)$ Suppose that there exist the framed curves $\mathcal{I}_1(\gamma,\nu_1,\nu_2), \mathcal{I}_2(\gamma,\nu_1,\nu_2):I \to \R^3 \times \Delta$, where $(\lambda^{\mathcal{I}}_1, \eta^{\mathcal{I}}_1)$, $(\lambda^{\mathcal{I}}_2, \eta^{\mathcal{I}}_2):I\to\R^2$ are smooth maps such that condition \eqref{involute_framed-condition} holds and $\theta^{\mathcal{I}}_1, \theta^{\mathcal{I}}_2:I\to\R$ are smooth functions. 
Moreover, suppose that $(\lambda^{\mathcal{I}}_1-\lambda^{\mathcal{I}}_2, \eta^{\mathcal{I}}_1-\eta^{\mathcal{I}}_2) \not\equiv (0, 0)$
Then the framed curves $\mathcal{I}_1(\gamma,\nu_1,\nu_2)$ and $\mathcal{I}_2(\gamma,\nu_1,\nu_2)$ are $(\nu_1,\nu_2,\overline{\mu})$-mates.
\end{lemma}
\demo
$(1)$ By Definition \ref{evolute_involute_framed} $(1)$, the framed curves $\mathcal{E}_1(\gamma,\nu_1,\nu_2), \mathcal{E}_2(\gamma,\nu_1,\nu_2):I \to \R^3 \times \Delta$ are given by
\begin{align*}
\mathcal{E}_1&(\gamma,\nu_1,\nu_2)(t)=(\Ev_1(\gamma,\nu_1,\nu_2)(t),\nu^{\mathcal{E}1}_1(\gamma,\nu_1,\nu_2)(t),\nu^{\mathcal{E}1}_2(\gamma,\nu_1,\nu_2)(t)),\\
&\Ev_1(\gamma,\nu_1,\nu_2)(t)=\gamma(t)+\lambda^{\mathcal{E}}_1(t)\nu_1(t)+\eta^{\mathcal{E}}_1(t)\nu_2(t),\\
&\nu^{\mathcal{E}1}_1(\gamma,\nu_1,\nu_2)(t)=\sin\theta^{\mathcal{E}}_1(t)\nu_1(t)+\cos\theta^{\mathcal{E}}_1(t)\nu_2(t),\\
&\nu^{\mathcal{E}1}_2(\gamma,\nu_1,\nu_2)(t)=\mu(t),\\
\mathcal{E}_2&(\gamma,\nu_1,\nu_2)(t)=(\Ev_2(\gamma,\nu_1,\nu_2)(t),\nu^{\mathcal{E}2}_1(\gamma,\nu_1,\nu_2)(t),\nu^{\mathcal{E}2}_2(\gamma,\nu_1,\nu_2)(t)),\\
&\Ev_2(\gamma,\nu_1,\nu_2)(t)=\gamma(t)+\lambda^{\mathcal{E}}_2(t)\nu_1(t)+\eta^{\mathcal{E}}_2(t)\nu_2(t),\\
&\nu^{\mathcal{E}2}_1(\gamma,\nu_1,\nu_2)(t)=\sin\theta^{\mathcal{E}}_2(t)\nu_1(t)+\cos\theta^{\mathcal{E}}_2(t)\nu_2(t),\\
&\nu^{\mathcal{E}2}_2(\gamma,\nu_1,\nu_2)(t)=\mu(t),
\end{align*}
where $\theta^{\mathcal{E}}_1, \theta^{\mathcal{E}}_2:I\to\R$ are smooth functions and $(\lambda^{\mathcal{E}}_1, \eta^{\mathcal{E}}_1), (\lambda^{\mathcal{E}}_2, \eta^{\mathcal{E}}_2):I\to\R^2$ are smooth maps such that condition \eqref{evolute_framed-condition} holds. 
If we take a smooth map $(\lambda, \eta):I\to\R^2$ by
\begin{align*}
\lambda(t)&=(-\lambda^{\mathcal{E}}_1(t)+\lambda^{\mathcal{E}}_2(t))\cos\theta^{\mathcal{E}}_1(t)-(-\eta^{\mathcal{E}}_1(t)+\eta^{\mathcal{E}}_2(t))\sin\theta^{\mathcal{E}}_1(t),\\
\eta(t)&=(-\lambda^{\mathcal{E}}_1(t)+\lambda^{\mathcal{E}}_2(t))\sin\theta^{\mathcal{E}}_1(t)+(-\eta^{\mathcal{E}}_1(t)+\eta^{\mathcal{E}}_2(t))\cos\theta^{\mathcal{E}}_1(t),
\end{align*}
then $(\lambda, \eta)\not\equiv(0,0)$. 
We denote $(\overline{\gamma},\overline{\nu}_1,\overline{\nu}_2)(t)=(\Ev_1(\gamma,\nu_1,\nu_2)(t)+\lambda(t)\mu^{\mathcal{E}1}(\gamma,\nu_1,\nu_2)(t)+\eta(t){\nu}^{\mathcal{E}1}_1(\gamma,\nu_1,\nu_2)(t),\sin\theta^{\mathcal{E}}_2\nu_1(t)+\cos\theta^{\mathcal{E}}_2\nu_2(t),\mu(t))$.
By a direct calculation, we have
\begin{align*}
\overline{\gamma}(t)&=\Ev_1(\gamma,\nu_1,\nu_2)(t)+\lambda(t)\mu^{\mathcal{E}1}(\gamma,\nu_1,\nu_2)(t)+\eta(t){\nu}^{\mathcal{E}1}_1(\gamma,\nu_1,\nu_2)(t)\\
&=\gamma(t)+\lambda^{\mathcal{E}}_2(t)\nu_1(t)+\eta^{\mathcal{E}}_2(t)\nu_2(t)\\
&=\Ev_2(\gamma,\nu_1,\nu_2)(t),
\end{align*}
Therefore, $(\overline{\gamma},\overline{\nu}_1,\overline{\nu}_2)$ is the framed curve and we have $\overline{\nu}_1(t)=\nu^{\mathcal{E}2}_1(\gamma,\nu_1,\nu_2)(t)$ and $\overline{\nu}_2(t)=\nu^{\mathcal{E}2}_2(\gamma,\nu_1,\nu_2)(t)$.
Therefore, the framed curves $(\overline{\gamma},\overline{\nu}_1,\overline{\nu}_2)=\mathcal{E}_2(\gamma,\nu_1,\nu_2)$ and $\mathcal{E}_1(\gamma,\nu_1,\nu_2)$ are $(\mu,\nu_1,\overline{\nu}_2)$-mates.
\par
$(2)$ We can also prove by the same method of $(1)$.
\enD
\par
We consider cases in which evolutes and involutes are not inverse operations.
By using Theorem \ref{evo_invo_framed_inverse} and Lemma \ref{relation_Invs_Evs}, we have the following.
\begin{theorem}\label{evo_invo_framed_inverse-another}
Let $(\gamma,\nu_1,\nu_2):I \to \R^3 \times \Delta$ be a framed curve with curvature $(\ell,m,n,\alpha)$.\par
$(1)$ Suppose that there exist the framed curves $\mathcal{I}(\gamma,\nu_1,\nu_2):I \to \R^3 \times \Delta$, where $(\lambda^{\mathcal{I}}, \eta^{\mathcal{I}}):I\to\R^2$ is a smooth map such that condition \eqref{involute_framed-condition} holds and $\theta^{\mathcal{I}}:I\to\R$ is a smooth function and $\mathcal{E}(\mathcal{I}(\gamma,\nu_1,\nu_2)):I \to \R^3 \times \Delta$, where $(\lambda^{\mathcal{E}}, \eta^{\mathcal{E}}):I \to \R^2$ with $(\lambda^{\mathcal{I}}+\lambda^{\mathcal{E}}\cos\theta^{\mathcal{I}}+\eta^{\mathcal{E}}\sin\theta^{\mathcal{I}}, \eta^{\mathcal{I}}-\lambda^{\mathcal{E}}\sin\theta^{\mathcal{I}}+\eta^{\mathcal{E}}\cos\theta^{\mathcal{I}})\not\equiv (0, 0)$ is a smooth map and $\theta^{\mathcal{E}}:I\to\R$ with $\theta^{\mathcal{E}}_1+\theta^{\mathcal{I}}\not\equiv0$ is a smooth function such that condition $(\ref{Ev(Inv)_framed})$ holds. 
Then the framed curves $(\gamma,\nu_1,\nu_2)$ and $\mathcal{E}(\mathcal{I}(\gamma,\nu_1,\nu_2))$ are $(\mu,\nu_1,\overline{\nu}_2)$-mates.
\par
$(2)$ Suppose that there exist the framed curves $\mathcal{E}(\gamma,\nu_1,\nu_2):I \to \R^3 \times \Delta$, where $(\lambda^{\mathcal{E}}, \eta^{\mathcal{E}}):I \to \R^2$ is a smooth map and $\theta^{\mathcal{E}}:I\to\R$ is a smooth function such that condition $(\ref{evolute_framed-condition})$ holds and $\mathcal{I}(\mathcal{E}(\gamma,\nu_1,\nu_2)):I \to \R^3 \times \Delta$, where $(\lambda^{\mathcal{I}}, \eta^{\mathcal{I}}):I \to \R^2$ with $(\lambda^{\mathcal{E}}+\lambda^{\mathcal{I}}\cos\theta^{\mathcal{E}}+\eta^{\mathcal{I}}\sin\theta^{\mathcal{E}}, \eta^{\mathcal{E}}-\lambda^{\mathcal{I}}\sin\theta^{\mathcal{E}}+\eta^{\mathcal{I}}\cos\theta^{\mathcal{E}})\not\equiv (0, 0)$ is a smooth map such that condition $(\ref{Inv(Ev)_framed})$ holds. 
Then the framed curves $(\gamma,\nu_1,\nu_2)$ and $\mathcal{I}(\mathcal{E}(\gamma,\nu_1,\nu_2))$ are $(\nu_1,\nu_2,\overline{\mu})$-mates.
\end{theorem}
\demo
$(1)$ By assumption, there exists the framed curve 
\begin{align*}
\mathcal{E}(\mathcal{I}(\gamma,\nu_1,\nu_2))(t)&=\left(\Ev(\mathcal{I}(\gamma,\nu_1,\nu_2)),\nu_1^{\mathcal{E}}(\mathcal{I}(\gamma,\nu_1,\nu_2)),\nu_2^{\mathcal{E}}(\mathcal{I}(\gamma,\nu_1,\nu_2))\right)(t),
\end{align*}
where $\theta^{\mathcal{E}}:I\to\R$ with $\theta^{\mathcal{E}}+\theta^{\mathcal{I}}\not\equiv0$ is a smooth function and $(\lambda^{\mathcal{E}}, \eta^{\mathcal{E}}):I\to\R^2$ with $(\lambda^{\mathcal{I}}+\lambda^{\mathcal{E}}\cos\theta^{\mathcal{I}}+\eta^{\mathcal{E}}\sin\theta^{\mathcal{I}}, \eta^{\mathcal{I}}-\lambda^{\mathcal{E}}\sin\theta^{\mathcal{I}}+\eta^{\mathcal{E}}\cos\theta^{\mathcal{I}})\not\equiv (0, 0)$ is a smooth map such that condition \eqref{Ev(Inv)_framed} holds. 
By Theorem \ref{evo_invo_framed_inverse} $(1)$, if we take $\lambda^{\mathcal{E}}_2(t)=\lambda^{\mathcal{I}}(t)\cos\theta^{\mathcal{I}}(t)-\eta^{\mathcal{I}}(t)\sin\theta^{\mathcal{I}}(t)$, $\eta^{\mathcal{E}}_2(t)=\lambda^{\mathcal{I}}(t)\sin\theta^{\mathcal{I}}(t)+\eta^{\mathcal{I}}(t)\cos\theta^{\mathcal{I}}(t)$ and $\theta^{\mathcal{E}}_2(t)=-\theta^{\mathcal{I}}(t)$, then condition (\ref{Ev(Inv)_framed}) also holds and we have $\mathcal{E}_2(\mathcal{I}(\gamma,\nu_1,\nu_2))(t)=(\gamma,\nu_1,\nu_2)(t)$.
By Lemma \ref{relation_Invs_Evs} $(1)$, the framed curves $\mathcal{E}(\mathcal{I}(\gamma,\nu_1,\nu_2))$ and $\mathcal{E}_2(\mathcal{I}(\gamma,\nu_1,\nu_2))=(\gamma,\nu_1,\nu_2)$ are $(\mu,\nu_1,\overline{\nu}_2)$-mates.
\par
$(2)$ We can also prove by the same method of $(1)$.
\enD

\begin{remark}\label{relation_theta-phi}{\rm
When $(\lambda^{\mathcal{I}}+\lambda^{\mathcal{E}}\cos\theta^{\mathcal{I}}+\eta^{\mathcal{E}}\sin\theta^{\mathcal{I}}, \eta^{\mathcal{I}}-\lambda^{\mathcal{E}}\sin\theta^{\mathcal{I}}+\eta^{\mathcal{E}}\cos\theta^{\mathcal{I}})\not\equiv (0, 0)$ (respectively, $(\lambda^{\mathcal{E}}+\lambda^{\mathcal{I}}\cos\theta^{\mathcal{E}}+\eta^{\mathcal{I}}\sin\theta^{\mathcal{E}}, \eta^{\mathcal{E}}-\lambda^{\mathcal{I}}\sin\theta^{\mathcal{E}}+\eta^{\mathcal{I}}\cos\theta^{\mathcal{E}})\not\equiv (0, 0)$) holds, then the set of regular points of $\mathcal{I}(\gamma,\nu_1,\nu_2)$ (respectively, $\mathcal{E}(\gamma,\nu_1,\nu_2)$), that is, $\alpha^{\mathcal{I}}\neq0$ (respectively, $\alpha^{\mathcal{E}}\neq0$) is a dense subset of $I$.
}
\end{remark}

We give a concrete example of the evolute and the involute of a framed curve.
\begin{example}{\rm
Let $(\gamma,\nu_1,\nu_2):[0,2\pi) \to \R^3 \times \Delta$ be 
\begin{align*}
\gamma(t)&=\left(\frac{28}{3}\cos^3 t,\frac{28}{3}\sin^3 t,-\frac{21}{4}\cos 2t\right), \\
\nu_1(t)&=(-\sin t,-\cos t,0), \ \nu_2(t)=\left(\frac{3}{5}\cos t,-\frac{3}{5}\sin t,\frac{4}{5}\right).
\end{align*}
By a direct calculation, $\mu(t)=\nu_1(t)\times\nu_2(t)= \left(-({4}\cos t)/{5},({4}\sin t)/{5},{3}/{5}\right)$
and $(\gamma,\nu_1,\nu_2)$ is a framed curve with the curvature $(\ell(t), m(t), n(t), \alpha(t))=(-{3}/{5}, {4}/{5}, 0, 35\cos t\sin t)$. 
If we take $\lambda^{\mathcal{I}}(t)=(125\cos 2t)/12$ and $\eta^{\mathcal{I}}(t)=(25\cos t\sin t)/3$, then condition \eqref{involute_framed-condition} is satisfies. Therefore, the framed curve $\mathcal{I}(\gamma,\nu_1,\nu_2)=(\mathcal{I}nv(\gamma,\nu_1,\nu_2), \nu^{\mathcal{I}}_1(\gamma,\nu_1,\nu_2),\nu^{\mathcal{I}}_2(\gamma,\nu_1,\nu_2))$ is given by
\begin{align*}
\mathcal{I}nv(\gamma,\nu_1,\nu_2)(t)&=\gamma(t)+\lambda^{\mathcal{I}}(t)\mu(t)+\eta^{\mathcal{I}}(t)\nu_1(t)=(\cos^3t,\sin^3t,\cos2t),\\
\nu^{\mathcal{I}}_1(\gamma,\nu_1,\nu_2)(t)&=\cos^{\mathcal{I}}(t)\mu(t)-\sin^{\mathcal{I}}(t)\nu_1(t), \\ 
\nu^{\mathcal{I}}_2(\gamma,\nu_1,\nu_2)(t)&=\sin^{\mathcal{I}}(t)\mu(t)+\cos^{\mathcal{I}}(t)\nu_1(t).
\end{align*}
Then $\mathcal{I}nv(\gamma,\nu_1,\nu_2)$ is an astroid. Moreover, if we take $\lambda^{\mathcal{E}}(t)=-({175}\cos t\sin t)/{4}\cos t\sin t$, $\eta^{\mathcal{E}}(t)=({875}\cos t\sin t)/{12}$ and $\theta^{\mathcal{E}}=\pi/2$, then condition \eqref{evolute_framed-condition} is satisfies. Therefore, the framed curve $\mathcal{E}(\gamma,\nu_1,\nu_2)=(\Ev(\gamma,\nu_1,\nu_2), \nu^{\mathcal{E}}_1(\gamma,\nu_1,\nu_2),\nu^{\mathcal{E}}_2(\gamma,\nu_1,\nu_2))$ is given by
\begin{align*}
&\Ev(\gamma,\nu_1,\nu_2)(t)=\gamma(t)+\lambda^{\mathcal{E}}(t)\mu(t)+\eta^{\mathcal{E}}(t)\nu_1(t)=\frac{637}{12}(\cos^3t,\sin^3t,\cos2t),\\
&\nu^{\mathcal{E}}_1(\gamma,\nu_1,\nu_2)(t)=\sin^{\mathcal{E}}(t)\nu_1(t)+\cos^{\mathcal{E}}(t)\nu_2(t)=\nu_1(t), \\
&\nu^{\mathcal{E}}_2(\gamma,\nu_1,\nu_2)(t)=\mu(t).
\end{align*}
Then $\Ev(\gamma,\nu_1,\nu_2)$ is also an astroid.
}
\end{example}


Nozomi Nakatsuyama, 
\\
Muroran Institute of Technology, Muroran 050-8585, Japan,
\\
E-mail address: 25096009b@muroran-it.ac.jp

\end{document}